\documentclass[12pt]{article}

\usepackage{latexsym,amsfonts,amsmath,amssymb,graphicx,hyperref}

\newcommand{\SA}{\mbox{\bf S}}
\newcommand{\HA}{\mbox{\bf H}}
\newcommand{\GA}{\mbox{\bf G}}
\newcommand{\TT}{\mbox{\bf T}}
\newcommand{\St}{\mathbb{S}}
\newcommand{\EE}{\mbox{\bf E}\,}
\newcommand{\PP}{\mbox{\bf P}\,}
\newcommand{\QQ}{\mbox{\bf Q}\,}

\newcommand{\R}{\mathbb{R}}
\newcommand{\C}{\mathbb{C}}

\newcommand{\N}{\mathbb{N}}

\newcommand{\Z}{\mathbb{Z}}
\newcommand{\Ree}{\mbox{Re}\,}
\newcommand{\Imm}{\mbox{Im}\,}
\newcommand{\A}{\mathbb{A}}
\newcommand{\CC}{\mbox{\bf C}}

\newcommand{\pa}{\partial}

\newcommand{\no}{\noindent}

\def\eps{\varepsilon}
\def\til{\widetilde}
\def\ha{\widehat}
\def\sem{\setminus}
\def\lin{\overline}
\def\vphi{\varphi}

\newtheorem{Lemma}{Lemma}[section]
\newtheorem{Theorem}{Theorem}[section]

\newtheorem{Proposition}{Proposition}[section]
\begin{document}
\title{\bf Some Properties of Annulus SLE}
\date{\today}
\author{Dapeng Zhan}
\maketitle
\begin{abstract} An annulus SLE$_\kappa$ trace
tends to a single point on the target circle, and the density
function of the end point satisfies some differential equation. Some
martingales or local martingales are found for annulus SLE$_4$,
SLE$_8$ and SLE$_{8/3}$. From the local martingale for annulus
SLE$_4$ we find a candidate of discrete lattice model that may have
annulus SLE$_4$ as its scaling limit. The local martingale for
annulus SLE$_{8/3}$ is similar to those for chordal and radial
SLE$_{8/3}$. But it seems that annulus SLE$_{8/3}$ does not satisfy
the restriction property.
\end{abstract}
%\vspace{-5mm}

\section{Introduction}
Schramm-Loewner evolution (SLE) is a family of random growth
processes  invented by O. Schramm in \cite{S-SLE} by connecting
Loewner differential equation with a one-dimentional Browinian
motion. SLE depend on a single parameter $\kappa\ge 0$, and behaves
differently for different value of $\kappa$. Schramm conjectured
that SLE$(2)$ is the  scaling limit of some loop-erased random walks
(LERW) and proved his conjuecture with some additional assumptions.
He also suggested that SLE$(6)$ and SLE$(8)$ should be the scaling
limits of certain discrete lattice models.

After Schramm's paper, there were many papers working on SLE.  In
the series of papers \cite{LSW1}\cite{LSW2}\cite{LSW3}, the locality
property of SLE$(6)$ was used to compute the intersection exponent
of plane Brwonian motion. In \cite{SS-6}, SLE$(6)$ was proved to be
the scaling limit of the cite percolation explorer on the triangle
lattice. It was proved in \cite{LSW-2} that SLE$(2)$ is the scaling
limit of the corresponding loop-erased random walk (LERW), and
SLE$(8)$ is the scaling limit of some uniform spanning tree (UST)
Peano curve. SLE$(4)$ was proved to be the scaling limit of the
harmonic exploer in \cite{SS}. SLE$(8/3)$ satisfies restriction
property, and was conjectured in \cite{LSW-8/3} to be the scaling
limit of some self avoiding walk (SAW). Chordal SLE$(\kappa,\rho)$
processes were also invented in \cite{LSW-8/3}, and they satisfy
one-sided restriction property. For basic properties of SLE, see
\cite{RS-basic}, \cite{LawSLE}, \cite{WerSLE}, \cite{WW-overall}.

The SLE invented by O. Schramm has a chordal and a radial version.
They are all defined in simply connected domains. In \cite{Zhan}, a
new version of SLE, called annulus SLE, was defined in doubly
connected domains as follows.

For $p>0$, let the annulus
$$\A_p=\{z\in\C:e^{-p}<|z|<1\},$$
and the circle
$$\CC_p=\{z\in\C:|z|=e^{-p}\}.$$ Then $\A_p$ is bounded by
$\CC_p$ and $\CC_0$. Let $\xi(t)$, $0\le t<p$, be a real valued
continuous function. For $z\in\A_p$, solve the annulus Loewner
differential equation
\begin{equation}
\partial_t\vphi_t(z)=\vphi_t(z)\SA_{p-t}(\vphi_t(z)/\exp(i\xi(t))),\,\,\,
0\le t<p,\,\,\,\vphi_0(z)=z,\label{annulus equation}
\end{equation}
where for $r>0$,
$$\SA_r(z)=\lim_{N\to\infty}\sum_{k=-N}^N\frac{e^{2kr}+z}{e^{2kr}-z}.$$
For $0\le t<p$, let $K_t$ be the set of $z\in\A_p$ such that the
solution $\vphi_s(z)$ blows up before or at time $t$. Then for each
$0\le t<p$, $\vphi_t$ maps $\A_p\sem K_t$ conformally onto
$\A_{p-t}$, and maps $\CC_p$ onto $\CC_{p-t}$. We call $K_t$
 and $\vphi_t$, respectively, $0\le t<p$, the annulus LE hulls and maps,
 respectively, of modulus $p$, driven by $\xi(t)$, $0\le t<p$. If
 $(\xi(t))=\sqrt{\kappa}B(t)$, $0\le t<p$, where $\kappa\ge 0$ and $B(t)$
is a standard linear Brownian motion, then $K_t$ and $\vphi_t$,
$0\le t<p$, are called standard annulus SLE$_\kappa$ hulls and maps,
respectively, of modulus $p$. Suppose $D$ is a doubly connected
domain with finite modulus $p$, $a$ is a boundary point and $C$ is a
boundary component of $D$ that does not contain $a$. Then there is
$f$ that maps $\A_p$ conformally onto $D$ such that $f(1)=a$ and
$f(\CC_p)=C$. Let $K_t$, $0\le t<p$, be standard annulus
SLE$_\kappa$ hulls. Then $(f(K_t),0\le t<p)$ is called an annulus
SLE$_\kappa(D;a\to C)$ chain.

It is known in \cite{Zhan} that annulus SLE$_\kappa$ is weakly
equivalent to radial SLE$_\kappa$. so from the existence of radial
SLE$_\kappa$ trace, we know the existence of a standard annulus
SLE$_\kappa$ trace, which is $\beta(t)=\vphi_t^{-1}(\exp(i\xi(t)))$,
$0\le t<p$. Almost surely $\beta$ is a continuous curve in
$\lin{\A_p}$, and for each $t\in[0,p)$, $K_t$ is the hull generated
by $\beta((0,t])$, i.e., the complement of the component of
$\A_p\sem\beta((0,t])$ whose boundary contains $\CC_p$. It is known
that when $\kappa=2$ or $\kappa=6$, $\lim_{t\to p}\beta(t)$ exists
and lies on $\CC_p$ almost surely. In this paper, we prove that this
is true for any $\kappa>0$. And we discuss the density function of
the distribution of the limit point. The density function should
satisfy some differential equation.

When $\kappa=2$, $8/3$, $4$, $6$, or $8$,  radial and chordal
SLE$_\kappa$ satisfy some special properties. Radial SLE$_6$
satisfies locality property. Since annulus SLE$_6$ is (strongly)
equivalent to radial SLE$_6$, so annulus SLE$_6$ also satisfies the
locality property. Annulus SLE$_2$ is the scaling limit of the
corresponding loop-erased random walk. In this paper, we discuss the
cases $\kappa=4$, $8$, and $8/3$. We find martingales or local
martingales for annulus SLE$_\kappa$ in each of these cases. From
the local martingale for annulus SLE$_4$, we may construct a
harmonic explorer whose scaling limit is annulus SLE$_4$. The
martingales for annulus SLE$_{8/3}$ are similar to the martingales
for radial and chordal SLE$_{8/3}$, which are used to show that
radial and chordal SLE$_{8/3}$ satisfy the restriction property.
However, the martingales for annulus SLE$_{8/3}$ does not help us to
prove that annulus SLE$_{8/3}$ satisfies  the restriction property.
On the contrary, it seems that annulus SLE$_{8/3}$ does not satisfy
the restriction property.

\section{Annulus Loewner Evolution in the Covering Space}
We often lift the annulus Loewner evolution to the covering space.
Let $e^i$ denote the map $z\mapsto e^{iz}$. For $p>0$, let
$\St_p=\{z\in\C:0<\Imm z<p\}$, $\R_p=ip+\R$, and
$\HA_p(z)=\frac1i\SA_p(e^i(z))$. Then $\St_p=(e^i)^{-1}(\A_p)$ and
$\R_p=(e^i)^{-1}(\CC_p)$. Solve
\begin{equation}\pa_t\til\vphi_t(z)=\HA_{p-t}(\til\vphi_t(z)-\xi(t)),
\,\,\,\til\vphi_0(z)=z. \label{annulus equation
covering}\end{equation}
 For $0\le t<p$, let $\til K_t$ be the set of
$z\in\St_p$ such that $\til\vphi_t(z)$ blows up before or at time
$t$. Then for each $0\le t<p$, $\til\vphi_t$ maps $\St_p\sem \til
K_t$ conformally onto $\St_{p-t}$, and maps $\R_p$ onto $\R_{p-t}$.
And for any $k\in\Z$, $\til\vphi_t(z+2k\pi)=\til\vphi_t(z)+2k\pi$.
We call $\til K_t$ and $\til \vphi_t$, $0\le t<p$, the annulus LE
hulls and maps, respectively, of modulus $p$ in the covering space,
driven by $\xi(t)$, $0\le t<p$. Then we have $\til
K_t=(e^i)^{-1}(K_t)$ and $e^i\circ\til\vphi_t=\vphi_t\circ e^i$. If
$(\xi(t))_{0\le t<p}$ has the law of $(\sqrt{\kappa}B(t))_{0\le
t<p}$, then $\til K_t$ and $\til \vphi_t$, $0\le t<p$, are called
standard annulus SLE$_\kappa$ hulls and maps, respectively, of
modulus $p$ in the covering space.

It is clear that $\HA_r$ is an odd function. It is analytic in $\C$
except at the set of simple poles $\{2k\pi+i2mr:k,m\in\Z\}$. And at
each pole $z_0$, the principle part is $\frac{2}{z-z_0}$. For each
$z\in\C$,  $\HA_r(z+2\pi)=\HA_r(z)$, and $\HA_r(z+i2r)=\HA_r(z)-2i$.
%Moreover, if $z\in\St_p$ and $\Imm z$ lies between $2k\pi$ and
% $(2k\pm 1)\pi$ for some $k\in\Z$, then $\pm\Ree \HA_r(z)>0$.

Let
$$f_r(z)=i\frac\pi r\HA_{\pi^2/r}(i\frac\pi r z).$$
Then $f_r$ is an odd function. It is analytic in $\C$ except at the
set of simple poles $\{z\in\C:i\frac\pi r z=2k\pi+i2m\pi^2/r\mbox{,
for some }k,m\in\Z\}=\{2m\pi-i2k r:m,k\in\Z\}$. And at each pole
$z_0$, the principle part is $\frac{2}{z-z_0}$. We then compute
$$f_r(z+2\pi)=i\frac\pi r\HA_{\pi^2/r}(i\frac\pi r z+i2\pi^2/r)=
i\frac\pi r(\HA_{\pi^2/r}(i\frac\pi r z)-2i)=f_r(z)+2\frac\pi r;$$
$$f_r(z+i2r)=i\frac\pi r\HA_{\pi^2/r}(i\frac\pi r z-2\pi)
=i\frac\pi r\HA_{\pi^2/r}(i\frac\pi r z)=f_r(z).$$ Let
$g_r(z)=f_r(z)-\HA_r(z)$. Then $g_r$ is an odd entire function, and
satisfies $$g_r(z+2\pi)=g_r(z)+2\pi/r,\quad g_r(z+i2r)=g_r(z)+2i,$$
for any $z\in\C$. Thus $g_r(z)=z/r$. So we have
\begin{equation}\HA_r(z)=f_r(z)-g_r(z)=i\frac\pi
r\HA_{\pi^2/r}(i\frac\pi r z)-\frac{z}r.\label{change HA}
\end{equation}

\section{Long Term Behaviors of Annulus SLE Trace}
In this section we fix $\kappa>0$ and $p>0$. Let $\vphi_t$ and
$K_t$, $0\le t<p$, be the annulus LE maps and hulls, respectively,
of modulus $p$ driven by $\xi(t)=\sqrt{\kappa}B(t)$, $0\le
 t<p$. Let $\til\vphi_t$ and $\til K_t$ be the corresponding annulus
LE maps and hulls in the covering space. Let $\beta(t)$ be the
corresponding annulus SLE$_\kappa$ trace.

 Let $Z_t(z)=\til\vphi_t(z)-\xi(t)$. Then we have
 $$dZ_t(z)=\HA_{p-t}(Z_t(z))dt-\sqrt{\kappa}dB(t).$$
 Let $W_t(z)=\frac{\pi}{p-t}Z_t(z)$. Then $W_t$ maps $(\St_p\sem \til K_t,\R_p)$
conformally onto $(\St_\pi,\R_\pi)$. From Ito's formula and equation
(\ref{change HA})
 we have
 $$dW_t(z)=\frac{\pi dZ_t(z)}{p-t}+\frac{\pi Z_t(z)}{(p-t)^2}=
-\frac{\pi
\sqrt{\kappa}dB(t)}{p-t}+\frac{\pi}{p-t}(\HA_{p-t}(Z_t(z))+\frac{Z_t(z)}{p-t})dt$$
$$=-\frac{\pi\sqrt{\kappa}dB(t)}{p-t}+\frac{\pi}{p-t}i\frac{\pi}{p-t}\HA_{\pi^2/(p-t)}
(i\frac{\pi}{p-t}Z_t(z))dt$$
$$=-\frac{\pi\sqrt{\kappa}dB(t)}{p-t}+\frac{i\pi^2}{(p-t)^2}\HA_{\pi^2/(p-t)}
(iW_t(z))dt.$$ Now we change variables as follows. Let
$s=u(t)=\pi^2/(p-t)$. Then $u'(t)=\pi^2/(p-t)^2$. For $\pi^2/p\le
s<\infty$, let $\ha W_s(z)=W_{u^{-1}(s)}(z)$. Then there is a
standard one dimensional Brownian motion $(B_1(s),s\ge \pi^2/p)$
such that
$$d\ha W_s(z)=\sqrt{\kappa}d B_1(s)+i\HA_s(i\ha W_s(z))ds,$$

Let $\ha\vphi_s(z)=\ha W_s(z)-\sqrt\kappa B_1(s)$. Then
$\pa_s\ha\vphi_s(z)=i\HA_s(i\ha W_s(z))$. Let $X_s(z)=\Ree\ha
W_s(z)$. For $z\in\R_p$, we have $\ha
W_s(z),\ha\vphi_s(z)\in\R_\pi$, so $\ha W_s(z)=X_s(z)+i\pi$. Thus
for $z\in\R_p$,
\begin{equation}\pa_s\Ree\ha\vphi_s(z)=\Ree\pa_s\ha\vphi_s(z)=\Ree(i\HA_s(i( X_s(z)+i\pi)))
=\lim_{M\to\infty}\sum_{k=-M}^M\frac{e^{X_s(z)}-e^{2ks}}{e^{X_s(z)}+e^{2ks}}
\label{equation on R_p}.\end{equation} Note that $\ha
W_s'(z)=\ha\vphi_s'(z)$. So for $z\in\R_p$,
\begin{equation*}\pa_s\ha\vphi_s'(z)=\sum_{k=-\infty}^\infty\frac{2e^{X_s(z)}e^{2ks}}
{(e^{X_s(z)}+e^{2ks})^2}\ha\vphi_s'(z),\end{equation*}
 which implies that
\begin{equation}\pa_s\ln|\ha\vphi_s'(z)|=\sum_{k=-\infty}^{+\infty}\frac{2e^{X_s(z)}e^{2ks}}
{(e^{X_s(z)}+e^{2ks})^2}.\label{logrithm
differentiation}\end{equation}

\begin{Lemma} For every $z\in\R_p$, $X_s(z)$ is not bounded on
$[\pi^2/p,\infty)$
 almost surely. \label{tend toinfinity}
\end{Lemma} {\bf Proof.} Suppose the lemma is not true. Then there is
$z_0\in\R_p$ and $a>0$ such that the probability that $|X_s(z_0)|<a$
for all $s\in[\pi^2/p,\infty)$ is positive. Let $X_s$ denote
$X_s(z_0)$. Then we have
$$dX_s=\sqrt{\kappa}d B_1(s)+\left(
\lim_{M\to\infty}\sum_{k=-M}^M\frac{e^{X_s}-e^{2ks}}{e^{
X_s}+e^{2ks}} \right)ds.$$ Let $T_a$ be the first time that
$|X_s|=a$. If such time does not exist, then let $T_a=\infty$. Let
$f(x)=\int_{-\infty}^x\cosh(s/2)^{-4/\kappa}ds$. Then $f$ maps $\R$
onto $(0,C(\kappa))$ for some $C(\kappa)<\infty$, and
 $f'(x)=\cosh(x/2)^{-4/k}$. So $f'(x)\frac{e^x-1}{e^x+1}+\frac\kappa
 2 f''(x)=0$. Let $U_s=f(X_s)$. Then
$$dU_s=f'(X_s)dX_s+\frac\kappa 2f''(X_s)ds$$$$=f'(X_s)\sqrt\kappa
d B_1(s)+f'(X_s)\lim_{M\to\infty}\left(
\sum_{k=-M}^1\frac{e^{X_s}-e^{2ks}}{e^{ X_s}+e^{2ks}}
+\sum_{k=1}^M\frac{e^{X_s}-e^{2ks}}{e^{ X_s}+e^{2ks}}\right)ds$$
$$=f'(X_s)\sqrt\kappa d B_1(s)+f'(X_s)\sum_{k=1}^\infty\frac{2\sinh(X_s)}
{\cosh(2ks)+\cosh(X_s)}ds.$$
 Let $v(s)=\int_{\pi^2/p}^s f'(X_t)^2dt$ for
$\pi^2/p\le s<T_a$. Let $\ha T_a=v(T_a)$. For $0\le r<\ha T_a$, let
$\ha U_r=U_{v^{-1}(r)}$. Then
$$d\ha U_r=\sqrt\kappa d B_2(r)+f'(X_{v^{-1}(r)})^{-1}\sum_{k=1}^\infty\frac{2\sinh(X_{v^{-1}(r)})}
{\cosh(2ks)+\cosh(X_{v^{-1}(r)})}dr,$$ where $ B_2(r)$ is  another
standard one dimensional Brownian motion. And $\ha T_a$ is a
stopping time w.r.t.\ $ B_2(r)$. Let
$$A(r)=f'(X_{v^{-1}(r)})^{-1}\sum_{k=1}^\infty\frac{2\sinh(X_{v^{-1}(r)})}
{\cosh(2ks)+\cosh(X_{v^{-1}(r)})};$$
$$M(r)=\exp\left(-\int_0^rA(s)\sqrt\kappa d B_2(s)-\frac\kappa
2\int_0^r A(s)^2ds\right).$$ For $0\le r<\ha T_a$,
$|X_{v^{-1}(r)}|<a$, so $|f'(X_{v^{-1}(r)})^{-1}|\le
\cosh(a/2)^{4/\kappa}$. And
$$\left|\sum_{k=1}^\infty\frac{2\sinh(X_{v^{-1}(r)})}
{\cosh(2ks)+\cosh(X_{v^{-1}(r)})}\right|\le\sum_{k=1}^\infty\frac{2\sinh(a)}
{e^{2ks}/2}=\frac{4\sinh(a)}{e^{2s}-1}.$$ Thus the Nivikov's
condition
$$\EE[\exp(\frac\kappa 2\int_0^{\ha T_a}A(s)^2ds)]<\infty$$
is satisfied. Let $\PP$ denote the original measure for $ B_2(r)$.
Define $\QQ$ on $\ha{\cal F}_{\ha T_a}$ such that
$d\QQ(\omega)=M_{\ha T_a}(\omega)d\PP(\omega)$. Then $(\ha U_r,0\le
r<\til T_a)$ is a one dimensional Brownian motion started from $0$
and stopped at time $\ha T_a$ w.r.t.\ the probability law $\QQ$. For
$0\le s<T_a$, $|X_s|\le a$, so $|f'(X_s|\ge \cosh(a/2)^{-4/\kappa}$.
Thus if $T_a=\infty$, then $\ha T_a=\infty$ too. From the hypothesis
of the proof, $\PP\{T_a=\infty\}>0$, so $\PP\{\ha T_a=\infty\}>0$.
Since $(\ha U_r,0\le r<\ha T_a)$ is a one dimensional Brownian
motion w.r.t.\ $\QQ$, so on the event that $\ha T_a=\infty$,
$\QQ\{\limsup_{r\to\infty}|\ha U_r|<\infty\}=0$. Thus
$\QQ\{\limsup_{r\to\infty}|\ha U_r|=\infty\}>0$. Since $\PP$ and
$\QQ$ are equivalent probability measures, so $\PP\{\limsup_{r\to
\ha T_a}|\ha U_r|=\infty\}>0$. Thus $\PP\{\limsup_{s\to
T_a}|U_s|=\infty\}>0$. This contradicts the fact that for all
$s\in[\pi^2/p,\infty)$, $U_s\in(0,C(\kappa))$ and
$C(\kappa)<\infty$. Thus the hypothesis is wrong, and the proof is
completed. $\Box$

\vskip 3mm

From this lemma and the definition of $X_t$, we know that for any
$z\in\R_p$, $(\Ree \til\vphi_t(z)-\sqrt\kappa B(t))/(p-t)$ is not
bounded on $t\in [0,p)$ a.s.. Since for any $k\in\Z$ and $z\in\R_p$,
$\til\vphi_t(z)-2k\pi=\til\vphi_t(z-2k\pi)$, so $(\Ree
\til\vphi_t(z)-2k\pi-\sqrt\kappa
B(t)))/(p-t)=(\Ree\til\vphi_t(z-2k\pi)-\sqrt\kappa B(t))/(p-t)$ is
not bounded on $t\in[0,p)$ a.s., which implies that $X_s(z)-2ks$ is
not bounded on $s\in[\pi^2/p,\infty)$ a.s..

\begin{Lemma} For every $z\in\R_p$, almost surely $\lim_{s\to\infty}
X_s(z)/s$ exists and the limit is an odd integer. \label{limit
exists}
\end{Lemma} {\bf Proof.} Fix $\eps_0\in(0,1/2)$ and $z_0\in\R_p$. Let
$X_s$ denote $X_s(z_0)$. There is $b>0$ such that the probability
that $|\sqrt\kappa B(t)|\le b+\eps_0t$ for any $t\ge 0$ is greater
than $1-\eps_0$. Since $\coth(x/2)\to\pm1$ as $x\in\R$ and
$x\to\pm\infty$, so there is $R>0$ such that when $\pm x\ge R$,
$\pm\coth(x/2)\ge 1-\eps_0$. Let $T=R+b+1$. If for any $s\ge 0$,
$|X_s-2ks|< T$ for some $k=k(s)\in\Z$, then there is $k_0\in\Z$ such
that $|X_s-2k_0 s|<T$ for all $s\ge T$. From the argument after
Lemma \ref{tend toinfinity}, the probability of this event is $0$.
Let $s_0$ be the first time that $|X_s-2ks|\ge T$ for all $k\in\Z$.
Then $s_0$ is finite almost surely. There is $k_0\in\Z$ such that
$2k_0s_0+T\le X_{s_0}\le 2(k_0+1)s_0-T$. Let $s_1$ be the first time
after $s_0$ such that $X_s=2k_0s+R$ or $X_s=2(k_0+1)s-R$. Let
$s_1=\infty$ if such time does not exist. For $s\in[s_0,s_1)$, we
have $X_s\in[2k_0s+R,2(k_0+1)s-R]$. Note that
$(e^x-e^{2ks})/(e^x+e^{2ks})\to\mp 1$ as $k\to\pm\infty$. So
$$\lim_{M\to\infty}\sum_{k=-M}^M\frac{e^{X_s}-e^{2ks}}{e^{X_s}+e^{2ks}}
=2k_0+\lim_{M\to\infty}\sum_{k=k_0-M}^{k_0+M}\frac{e^{X_s}-e^{2ks}}{e^{X_s}+e^{2ks}}$$
$$=2k_0+\lim_{M\to\infty}\sum_{j=-M}^{M}\frac{e^{X_s-2k_0s}-e^{2js}}
{e^{X_s-2k_0s}+e^{2js}}$$$$=2k_0+\coth(\frac{X_s-2k_0s}2)+\sum_{j=1}^\infty
\frac{2\sinh(X_s-2k_0s)}{\cosh(2js)+\cosh(X_s-2k_0s)}$$
$$\ge 2k_0+\coth(\frac{X_s-2k_0s}2)\ge 2k_0+1-\eps_0;$$
and
$$\lim_{M\to\infty}\sum_{k=-M}^M\frac{e^{X_s}-e^{2ks}}{e^{X_s}+e^{2ks}}
=2(k_0+1)+\lim_{M\to\infty}\sum_{j=-M}^{M}\frac{e^{X_s-2(k_0+1)s}-e^{2js}}
{e^{X_s-2(k_0+1)s}+e^{2js}}$$
$$=2k_0+2+\coth(\frac{X_s-2(k_0+1)s}2)+\sum_{j=1}^\infty
\frac{2\sinh(X_s-2(k_0+1)s)}{\cosh(2js)+\cosh(X_s-2(k_0+1)s)}$$
$$\le 2k_0+2+\coth(\frac{X_s-2(k_0+1)s}2)\le
2k_0+2+(-1+\eps_0)=2k_0+1+\eps_0.$$ From equation (\ref{equation on
R_p}), we have that for  $s\in[s_0,s_1)$,
$$(2k_0+1-\eps_0)(s-s_0)\le \Ree\ha\vphi_s(z_0)-\Ree\ha\vphi_{s_0}(z_0)
\le(2k_0+1+\eps_0)(s-s_0).$$ Note that
$X_s=\Ree\ha\vphi_s(z_0)-\sqrt\kappa B_1(s)$, and $(\sqrt\kappa
B_1(s)-\sqrt\kappa B_1(s_0),s\ge s_0)$ has the same distribution as
$(\sqrt\kappa B(s-s_0),s\ge s_0)$. Let $E_b$ denote the event that
$|\sqrt\kappa B_1(s)-\sqrt\kappa B_1(s_0)|\le b+\eps_0(s-s_0)$ for
all $s\ge s_0$. Then $\PP(E)>1-\eps_0$. And on the event $E_b$, we
have
$$(2k_0+1-\eps_0)(s-s_0)-b-\eps_0(s-s_0)\le X_s-X_{s_0}$$$$\le(2k_0+1+\eps_0)
(s-s_0)+b+\eps_0(s-s_0),$$ from which follows that
$$ X_s\le X_{s_0}+(2k_0+1+\eps_0)(s-s_0)+b+\eps_0(s-s_0)$$$$
\le 2(k_0+1)s_0-T+(2k_0+1+\eps_0)(s-s_0)+b+\eps_0(s-s_0)$$
$$= 2(k_0+1)s-T+b-(1-2\eps_0)(s-s_0)\le 2(k_0+1)s-R-1$$
and $$X_s\ge X_{s_0}+(2k_0+1-\eps_0)(s-s_0)-b-\eps_0(s-s_0)$$
$$\ge 2k_0s_0+T+(2k_0+1-\eps_0)(s-s_0)-b-\eps_0(s-s_0)$$
$$=2k_0s+T-b+(1-2\eps_0)(s-s_0)\ge 2k_0s +R+1.$$
So  on the event $E_b$ we have $s_1=\infty$, which implies that
$2k_0s+R\le X_s\le 2(k_0+1)s-R$  for all $s\ge s_0$, and so
$\pa_s\Ree\ha\vphi_s(z_0)\in(2k_0+1-\eps_0,2k_0+1+\eps_0)$ for all
$s\ge s_0$. Thus the event that
$$2k_0+1-\eps_0\le\liminf_{s\to\infty}\Ree\ha\vphi_s(z_0)/s
\le\limsup_{s\to\infty}\Ree\ha\vphi_s(z_0)/s\le 2k_0+1+\eps_0$$ has
probability greater than $1-\eps_0$. Since we may choose $\eps_0>0$
arbitrarily small, so a.s.\ $\lim_{s\to\infty}\Ree\ha\vphi_s(z_0)/s$
exists and the limit is $2k_0+1$ for some $k_0\in\Z$. The proof is
now finished by the facts that
$X_s(z_0)=\Ree\ha\vphi_s(z_0)+\sqrt\kappa B_1(s)$ and
$\lim_{s\to\infty} B_1(s)/s=0$. $\Box$

\vskip 3mm

Let $$m_-=\sup\{x\in\R:\lim_{s\to\infty}X_s(x+ip)/s\le -1\}$$ and
$$m_+=\inf\{x\in\R:\lim_{s\to\infty}X_s(x+ip)/s\ge 1\}.$$
Since $X_s(x_1+ip)<X_s(x_2+ip)$ if $x_1<x_2$, so we have $m_-\le
m_+$. If the event that $m_-<m_+$ has a positive probability, then
there is $a\in\R$ such that the event that $m_-<a<m_+$ has a
positive probability. From the definitions, $m_-<a<m_+$ implies that
$\lim_{s\to\infty}X_s(a+ip)/s\in(-1,1)$, which is an event with
probability $0$ by Lemma \ref{limit exists}. This contradiction
shows that $m_-=m_+$ a.s.. Let $m=m_+$. For any $t\in[0,p)$,
$z\in\St_p\sem\til K_t$ and $k\in\Z$, since
$\til\vphi_t(z+2k\pi)=\til\vphi_t(z)+2k\pi$, so
$Z_t(z+2k\pi)=Z_t(z)+2k\pi$, then we have
$W_t(z+2k\pi)=W_t(z)+2k\pi^2/(p-t)$. Thus $X_s(z+2k\pi)=X_s(z)+2ks$
for any $s\in[\pi^2/p,\infty)$, $z\in\St_p\sem\til K_{p-\pi^2/s}$
and $k\in\Z$. If $x\in(m+2k\pi,m+2(k+1)\pi)$ for some $k\in\Z$, then
$x-2k\pi>m$ and $x-2(k+1)\pi<m$. So
$$\lim_{s\to\infty}X_s(x+ip)/s=\lim_{s\to\infty}(X_s(x-2k\pi+ip)+2ks)/s
\ge 2k+1$$ and
$$\lim_{s\to\infty}X_s(x+ip)/s=\lim_{s\to\infty}(X_s(x-2(k+1)\pi+ip)+2(k+1)s)/s
\le 2k+1.$$ Therefore $\lim_{s\to\infty}X_s(x+ip)/s=2k+1$.

Let $K_p=\cup_{0\le t<p} K_t$ and $\til K_p=\cup_{0\le t<p}\til
K_t$. Then $K_p=e^i(\til K_p)$, and so $\lin{K_p}=e^i(\lin{\til
K_p})$.

\begin{Lemma} $\lin{K_p}\cap\CC_p=\{e^{-p+im}\}$ almost
surely. \label{Intersection tilde}
\end{Lemma}
{\bf Proof.} We first show that $m+ip\in\lin{\til K_p}$. If this is
not true, then there is $a,b>0$ such that the distance between
$[m-a+ip,m+a+ip]$ and $\til K_t$ is greater than $b$ for all
$t\in[0,p)$. From the definition of $m$, we have $X_s(m\pm
a+ip)\to\pm\infty$ as $s\to\infty$. Thus $\Ree \ha\vphi_s(m+
a+ip)-\Ree\ha\vphi_s(m-a+ip)\to \infty$ as $s\to\infty$. So there is
$c(s)\in(m-a,m+a)$ such that $\ha\vphi_s'(c(s)+ip)\to\infty$ as
$s\to\infty$. Since $\ha\vphi_s$ maps $(\St_p\sem\til
K_{p-\pi^2/s},\R_p)$ conformally onto $(\St_\pi,\R_\pi)$, so by
Koebe's $1/4$ theorem, the distance between $c(s)+ip$ and $\til
K_{p-\pi^2/s}$ tends to $0$ as $s\to\infty$. This is a
contradiction. Thus $m+ip\in\lin{\til K_p}$.

Now fix $x_1<x_2\in(m,m+2\pi)$. Then $X_s(x_j+ip)/s\to 1$ as
$s\to\infty$ for $j=1,2$. So there is $s_0$ such that
$X_s(x_j+ip)\in(s/2,3s/2)$ for $s\ge s_0$ and $j=1,2$. So if
$x_0\in[x_1,x_2]$ and $s\ge s_0$, then $X_s(x_0+ip)\in(s/2,3s/2)$,
and so
$$\sum_{k=-\infty}^{+\infty}\frac{e^{X_s(x_0+ip)}e^{2ks}}
{(e^{X_s(x_0+ip)}+e^{2ks})^2}\le\sum_{k=-\infty}^{0}e^{2ks-X_s(x_0+ip)}+
\sum_{k=1}^{+\infty} e^{X_s(x_0+ip)-2ks}$$
$$\le\sum_{k=-\infty}^{0}e^{2ks-s/2}+
\sum_{k=1}^{+\infty}
e^{3s/2-2ks}=\frac{2e^{-s/2}}{1-e^{-2s}}\le\frac{2e^{-s/2}}
{1-e^{-2\pi^2/p}}.$$ From equation (\ref{logrithm differentiation}),
for all $s\ge s_0$,
$$\pa_s\ln|\ha\vphi_s'(x_0+ip)|\le\frac{4e^{-s/2}}
{1-e^{-2\pi^2/p}},$$ which implies that
$$\ln|\ha\vphi_s'(x_0+ip)|\le \ln|\ha\vphi_{s_0}'(x_0+ip)|+\frac{8e^{-s_0/2}}
{1-e^{-2\pi^2/p}}.$$ So there is $M<\infty$ such that
$|\ha\vphi_s'(x_0+ip)|\le M$ for all $x_0\in[x_1,x_2]$ and $s\ge
s_0$. From Koebe's $1/4$ theorem, we see that $\til K_{t}$ is
uniformly bounded away from $[x_1+ip,x_2+ip]$ for $t\in[0,p)$. Thus
$[x_1+ip,x_2+ip]\cap\lin{\til K_p}=\emptyset$. Since $x_1<x_2$ are
chosen arbitrarily from $(m,m+2\pi)$, so $(m+ip,m+2\pi+ip)\cap
\lin{\til K_p}=\emptyset$. Thus $\lin{\til K_p}\cap
[m+ip,m+2\pi+ip)=\{m+ip\}$. Since $\CC_p=e^i([m+ip,m+2\pi+ip))$, so
$\lin{K_p}\cap\CC_p=\{e^i(m+ip)\}=\{e^{-p+im}\}$. $\Box$

\begin{Lemma} For every $\eps\in(0,1)$, there is $C_0>0$
depending on $\eps$
 such that if $q\in(0,\frac{2\pi^2}{\ln(2)}]$, and $L_t$, $0\le t<q$, are
  standard annulus SLE$_\kappa$ hulls of modulus $q$, then the
probability that $\cup_{0\le t<q} L_t\subset\{e^{iz}:|\Ree z|\le
C_0q\}$ is greater than $1-\eps$. \label{bounded probability}
\end{Lemma}
{\bf Proof.} Let $q_0=\frac{2\pi^2}{\ln(2)}$. Suppose $q\in(0,q_0]$.
Let $L_t$ and $\psi_t$, $0\le t<q$, be the annulus LE hulls and maps
of modulus $q$ driven by $\sqrt\kappa B(t)$, $0\le t<q$. Let $\til
L_t$ and $\til\psi_t$, $0\le t<q$,  be the corresponding annulus LE
hulls and maps in the covering space. There is $b=b(\eps)>0$ such
that the probability that $|\sqrt\kappa B(t)|\le b+t/4$ for all
$t\ge 0$ is greater than $1-\eps$. Let $R=\ln(64)$ and
$C_0=(R+b+1)/\pi$. Let $s_0=\pi^2/q$. Let
$Z_t(z)=\til\psi_t(z)-\sqrt\kappa B(t)$, $W_t(z)=\pi Z_t(z)/(q-t)$
for $0\le t<q$. Let $\ha W_s(z)=W_{q-\pi^2/s}(z)$ for $s_0\le
s<\infty$. Then there is another standard one dimensional Brownian
motion $B_1(s)$, $s\ge s_0$,  such that $\ha\psi_s$ defined by
$\ha\psi_s(z)=\ha W_s(z)+\sqrt\kappa B_1(s)$ satisfies
$$\pa_s\ha\psi_s(z)=\lim_{M\to\infty}\sum_{k=-M}^M\frac{e^{\ha W_s(z)}
+e^{2ks}}{e^{\ha W_s(z)}-e^{2ks}}$$
 for $s_0\le s<\infty$. Let $E_\eps$ be the event that $|\sqrt\kappa
B_1(s)-\sqrt\kappa B_1(s_0)|\le b+(s-s_0)/4$ for all $s\ge s_0$.
Then $\PP(E_\eps)>1-\eps$. Fix $z_0\in\St_q$ with $C_0q<\Ree
z_0<2\pi-C_0q$. We claim that in the event $E_\eps$,
$\til\psi_t(z_0)$ never blows up for $0\le t<q$. If this claim is
justified, then on the event $E_\eps$, $z_0\not\in\til L_t$ for any
$0\le t<q$ and $z_0\in\St_q$ with $C_0q<\Ree z_0<2\pi-C_0q$. So
$\cup_{0\le t<q}\til L_t$ is disjoint from $\{z\in\C:C_0q<\Ree z<
2\pi-C_0q\}$. Since $L_t=e^i(\til L_t)$, so $\cup_{0\le t<q}L_q$ is
disjoint from $\{e^{iz}:C_0q<\Ree z< 2\pi-C_0q\}$ on the event
$E_\eps$. Then we are done.

Assume the event $E_\eps$. Let $Z_t$ denote $Z_t(z_0)$, $W_t$ denote
$W_t(z_0)$, $\ha W_s$ denote $\ha W_s(z_0)$, and $\ha\psi_s$ denote
$\ha\psi_s(z_0)$. If $\til\psi_t(z_0)$ blows up at time $t_*<q$,
then $Z_t\to 2k\pi$ for some $k\in\Z$ as $t\to t_*$. Then $\ha W_s-
2ks\to 0$ as $s\to \pi^2/(q-t_*)$.  Since $\Ree Z_0=\Ree
z_0\in[C_0q,2\pi-C_0q]$, so $\ha
W_{s_0}=W_0\in[C_0\pi,2s_0-C_0\pi]\subset(R,2s_0-R)$, and so there
is a first time $s_1>s_0$ such that  $\Ree\ha
W_{s_1}\in\{R,2s_1-R\}$. Then for $s\in[s_0,s_1]$, we have $\Ree\ha
W_{s}\in[R,2s-R]$.
 Then
$$|\lim_{M\to\infty}\sum_{k=-M}^M\frac{e^{\ha W_s}+e^{2ks}}
{e^{\ha W_s}-e^{2ks}}-1|\le\sum_{k=-\infty}^{0}|\frac{e^{\ha
W_s}+e^{2ks}} {e^{\ha
W_s}-e^{2ks}}-1|+\sum_{k=1}^\infty|\frac{e^{\ha W_s}+e^{2ks}}
{e^{\ha W_s}-e^{2ks}}+1|$$$$\le\sum_{k=-\infty}^{0}\frac{2} {|e^{\ha
W_s-2ks}|-1}+\sum_{k=1}^\infty\frac{2} {|e^{2ks-\ha
W_s}|-1}\le\sum_{k=-\infty}^{0}\frac
4{e^{R-2ks}}+\sum_{k=1}^\infty\frac{4} {e^{2ks-(2s-R)}}$$$$\le
\frac{8e^{-R}}{1-e^{-2s}}\le 16e^{-R}\le \frac 14,$$ where we use
the fact that $e^{-R}\le \frac 1{64}$ and $e^{-2s}\le
e^{-2s_0}=e^{-2\pi^2/q}\le e^{-2\pi^2/q_0}\le\frac 12$. Thus
$$|(\ha W_{s_1}-\ha W_{s_0})-(s_1-s_0)|\le
|(\ha\psi_{s_1}-\ha\psi_{s_0})-(s_1-s_0)|+|\sqrt\kappa
B_1(s_1)-\sqrt\kappa B_1({s_0})|$$
$$\le(s_1-s_0)/4+b+(s_1-s_0)/4=b+(s_1-s_0)/2.$$
Then we have
$$\Ree\ha W_{s_1}\ge\Ree\ha W_{s_0}+(s_1-s_0)-b-(s_1-s_0)/2
\ge C_0\pi+(s_1-s_0)/2-b>R$$ and $$\Ree\ha W_{s_1}\le \Ree\ha
W_{s_0}+(s_1-s_0)+b+(s_1-s_0)/2\le
2s_0-C_0\pi+b+3(s_1-s_0)/2$$$$=2s_1-(s_1-s_0)/2-C_0\pi+b<2s_1-R.$$
This contradicts that $\Ree\ha W_{s_1}\in\{R,2s_1-R\}$. Thus
$\til\psi_t(z_0)$ does not blow up for $t\in [0,q)$. Then the claim
is justified, and the proof is finished. $\Box$

\vskip 3mm

For two nonempty sets $A_1,A_2\subset\A_p$, we define the angular
distance between $A_1$ and $A_2$ to be $d_a(A_1,A_2)=\inf\{|\Ree
z_1-\Ree z_2|:e^{iz_1}\in A_1,e^{iz_2}\in A_2\}$. For a nonempty set
$A\subset\A_p$, we define the angular diameter of $A$ to be
$diam_a(A)=\sup\{d_a(z_1,z_2):z_1,z_2\in A\}$. If $A$ intersects
both $A_1$ and $A_2$, then $d_a(A_1,A_2)\le diam_a(A)$. In the above
lemma, $\cup_{0\le t<q}L_t\subset\{e^{iz}:|\Ree z|\le C_0q\}$
implies that $diam_a(\cup_{0\le t<q}L_t)\le 2C_0q$. Form conformal
invariance  and comparison principle of extremal distance, we have
that for any $d>0$, there is $h(d)>0$ such that for any $p>0$, if
for $j=1,2$, $A_j$ is a union of connected subsets of $\A_p$, each
of which touches both $\CC_p$ and $\CC_0$, and the extremal distance
between $A_1$ and $A_2$ in $\A_p$ is greater than $h(d)$, then
$d_a(A_1,A_2)>dp$.

\begin{Theorem} $\lim_{t\to p}\beta(t)=e^{-p+im}$ almost
surely.\label{limit}
\end{Theorem}
{\bf Proof.} From Lemma \ref{Intersection tilde}, the distance from
$e^{-p+im}$ to $K_t$ tends to $0$ as $t\to p$ a.s.. Since $K_t$ is
the hull generated by $\beta((0,t])$, so  the distance from
$e^{-p+im}$ to $\beta((0,t])$ tends to $0$ as $t\to p$ a.s.. Suppose
the theorem does not hold. Then there is $a,\delta>0$ such that the
event that $\limsup_{t\to p}|e^{-p+im}-\beta(t)|>a$ has probability
greater than $\delta$. Let $E_1$ denote this event. Let
$\eps=\delta/4$. Let $C_0$ depending on $\eps$ be as in Lemma
\ref{bounded probability}. Let $R=\min\{a,e^{-p}\}$ and
$r=\min\{1-e^{-p},R\exp(-2\pi h(2C_0+1))\}$, where $h$ is the
function in the argument before this theorem. Since $K_t$ is
generated by $\beta((0,t])$, and $e^{-p+im}\in\lin{K_p}$ a.s., so
the distance between $e^{-p+im}$ and $\beta((0,t])$ tends to $0$
a.s.\ as $t\to p$. So there is $t_0\in(0,p)$ such that the event
that the distance between $e^{-p+im}$ and $\beta((0,t_0])$ is less
than $r$ has probability greater than $1-\eps$. Let $E_2$ denote
this event. Let $q_0=\frac{2\pi^2}{\ln(2)}$,
$T=\max\{t_0,p-q_0,-\ln(r+e^{-p})\}$, $p_T=p-T$, and
$\xi_T(t)=\xi(T+t)-\xi(T)$ for $0\le t<p_T$. Let
$K_{T,t}=\vphi_{T}(K_{T+t}\sem K_{T})/e^{i\xi(T)}$ and
$\vphi_{T,t}(z)=\vphi_{T+t}\circ\vphi_T^{-1}(\exp(i\xi(T))z)/\exp(i\xi(T))$
for $0\le t<p_T$. Then one may check that $K_{T,t}$ and
$\vphi_{T,t}$, $0\le t<p_T$, are the annulus LE hulls and maps of
modulus $p_T$ driven by $\xi_T$. Since $\xi_T(t)$ has the same law
as $\sqrt\kappa B(t)$ and $p_T=p-T\le q_0$, so from Lemma
\ref{bounded probability}, the event that $diam_a(\cup_{0\le
t<p_T}K_{T,t})\le 2C_0p_T$ has probability greater than $1-\eps$.
Let $E_3$ denote this event. Since
$\PP(E_1^c)+\PP(E_2^c)+\PP(E_3^c)<(1-\delta)+\eps+\eps<1$, so
$\PP(E_1\cap E_2\cap E_3)>0$. This means that the events $E_1$,
$E_2$ and $E_3$ can happen at the same time. We will prove that this
is a contradiction. Then the theorem is proved.

Assume the event $E_1\cap E_2\cap E_3$. Let $A_r$  ($A_R$, resp.) be
the union of connected components of
$\{z\in\C:|z-e^{-p+im}|=r\}\cap(\A_p\sem K_{T})$
($\{z\in\C:|z-e^{-p+im}|=R\}\cap(\A_p\sem K_{T})$, resp.) that touch
$\CC_p$. From the properties of $\beta$ in the event $E_1$ and
$E_2$, we see that $A_r$ and $A_R$ both intersect $K_p\sem K_{T}$.
Since the distance between $e^{-p+im}$ and $K_{T}$ is less than $r$,
and $r<R$, so both $A_r$ and $A_R$ are unions of two curves which
touch both $\CC_p$ and $\CC_0\cup K_{T}$. Let
$B_r=e^{-i\xi(T)}\vphi_{T}(A_r)$ and
$B_R=e^{-i\xi(T)}\vphi_{T}(A_R)$. Then both $B_r$ and $B_R$ are
unions of two curves in $\A_{p_T}$ that touch both $\CC_{p_T}$ and
$\CC_0$.

The extremal distance between $A_r$ and $A_R$ in $\A_p\sem K_{T}$ is
at least $\ln(R/r)/(2\pi)\ge h(2C_0+1)$. Thus the extremal distance
between $B_r$ and $B_R$ in $\A_{p_T}$ is at least $h(2C_0+1)$. So
the angular distance between $B_r$ and $B_R$ is at least
$(2C_0+1)p_T$. Since $A_R$ and $A_r$ both intersect $K_p\sem K_{T}$,
so $B_R$ and $B_r$ both intersect $\vphi_{T}(K_p\sem
K_{T})/e^{i\xi(T)}=\cup_{0\le t<p_T}K_{T,t}$, which implies that
$diam_a(\cup_{0\le t<p_T}K_{T,t})\ge (2C_0+1)p_T$. However, in the
event $E_3$, $diam_a(\cup_{0\le t<p_T}K_{T,t})\le 2C_0p_T$. This
contradiction finishes the proof. $\Box$

\vskip 3mm

Let's see what can we say about the distribution of $\lim_{t\to
p}\beta(t)$. Let $\til\beta(t)=\til\vphi_t^{-1}(\xi(t))$. Then
$\til\beta$ is a simple curve in $\St_p$ started from $0$, and
$\beta(t)=e^i(\til\beta(t))$. From Theorem \ref{limit}, $\lim_{t\to
p}\til\beta(t)$ exists and lies on $\R_p$. We call $\til\beta$ an
annulus SLE$_\kappa$ trace in the covering space. Let $m_p+ip$
denote the limit point, where $m_p$ is a real valued random
variable.

Suppose the distribution of $m_p$ is absolutely continuous w.r.t.\
the Lebesgue measure, and the density function $\til\lambda(p,x)$ is
$C^{1,2}$ continuous. This hypothesis is very likely to be true, but
the proof is still missing now. We then have
$\int_\R\til\lambda(p,x)dx=1$ for any $p>0$. Since the distribution
of $\til\beta$ is symmetric w.r.t.\ the imaginary axis, so is the
distribution of $\lim_{t\to p}\til\beta(t)$. Thus
$\til\lambda(p,-x)=\til\lambda(p,x)$. Moreover, we expect that when
$p\to 0$ the distribution of $(m_p+ip)*\frac\pi p$ tends to the
distribution of the limit point of a strip SLE$_\kappa$ trace
introduced in \cite{thesis}, whose density is
$\cosh(x/2)^{-4/\kappa}/C(\kappa)$ for some $C(\kappa)>0$. If this
is true, then the distribution of $m_p$ tends to the point mass at
$0$ as $p\to 0$.

For $0\le t<p$, let ${\cal F}_t$ be the $\sigma-$algebra generated
by $\xi(s)$, $0\le s\le t$. Fix $T\in[0,p)$. Let $p_T=p-T$. For
$0\le t<p_T$, let $\xi_T(t)=\xi(T+t)-\xi(T)$. Then $\xi_T(t)$ has
the same distribution as $\sqrt\kappa B(t)$, and is independent of
${\cal F}_T$. For $0\le t<T$, let
$$\til\vphi_{T,t}(z)=\til\vphi_{T+t}\circ\til\vphi_T^{-1}(z+\xi(T))-\xi(T).$$
Then
$\pa_t\til\vphi_{T,t}(z)=\HA_{p_T-t}(\til\vphi_{T,}(z)-\xi_T(t))$,
and $\til\vphi_{T,0}(z)=z$. Thus $\til\vphi_{T,t}(z)$, $0\le t<p_T$,
are annulus LE maps of modulus $p_T$ in the covering space driven by
$\xi_T(t)$, $0\le t<p_T$, and so are independent of ${\cal F}_T$.
Let
\begin{equation}\til\beta_T(t)=\til\vphi_{T,t}^{-1}(\xi_T(t))=\til\vphi_T\circ
\til\vphi_{T+t}^{-1}(\xi(T_t))-\xi(T)=\til\vphi_T(\til\beta(T+t))-\xi(T),
\label{trace change}\end{equation} for $0\le t<p_T$. Then
$\til\beta_T(t)$, $0\le t<p_T$, is a standard annulus SLE$_\kappa$
trace of modulus $p_T$ in the covering space, and is independent of
${\cal F}_T$. Thus $\lim_{t\to p_T}\beta_T(t)$ exists and lies on
$\R_{p_T}$ a.s.. Let $m_{p_T}+ip_T$ denote the limit point. Then
$m_{p_T}$ is independent of ${\cal F}_T$, and the density of
$m_{p_T}$ w.r.t.\ the Lebesgue measure is $\til\lambda(p_T,\cdot)$.
From equation (\ref{trace change}), we see
$m_{p_T}=\til\vphi_T(m_p+ip)-ip_T-\xi(T)$. For $0\le t<p$, let
$\til\psi_t(z)=\til\vphi_t(z+ip)-i(p-t)$. Then $\til\psi_t$ takes
real values on $\R$, and
$\pa_t\til\psi_t(z)=\ha\HA_{p-t}(\til\psi_t(z)-\xi(t))$, where
$\ha\HA_r(z)=\HA(z+ir)+i$. Let $X_t(z)=\til\psi_t(z)-\xi(t)$ for
$0\le t<p_T$. So $m_{p_T}=X_T(m_p)$. From the differential equation
for $\til\psi_t$, we get
$$dX_t(x)=\ha\HA_{p-t}(X_t(x))dt-d\xi(t);$$
and $$dX_t'(x)=\ha\HA_{p-t}'(X_t(x))X_t'(x)dt.$$

Let $a<b\in\R$. Then
$\{m_p\in[a,b]\}=\{m_{p_T}\in[X_T(a),X_T(b)]\}$. Since $m_{p_T}$ has
density $\til\lambda(p_T,\cdot)$ and is independent of ${\cal F}_T$,
and $X_T$ is ${\cal F}_T$ measurable, so
$$\EE[{\bf 1}_{\{m_p\in[a,b]\}}|{\cal F}_T]=\int_{X_T(a)}^{X_T(b)}
\til\lambda(p-T,x)dx=\int_a^b\til\lambda(p-T,X_T(x))X_T'(x)dx.$$
Thus $(\int_a^b\til\lambda(p-t,X_t(x))X_t'(x)dx, 0\le t<p)$ is a
martingale w.r.t.\ $\{{\cal F}_t\}_{t=0}^p$. Fix $x\in\R$. Choose
$a<x<b$ and let $a,b\to x$. Then
$(\til\lambda(p-t,X_t(x))X_t'(x),0\le t<p)$ is a martingale w.r.t.\
$\{{\cal F}_t\}_{t=0}^p$. From Ito's formula, we have
\begin{equation}-\pa_1\til\lambda(r,x)+\ha\HA_r'(x)\til\lambda(r,x)+
\ha\HA_r(x)\pa_2\til\lambda(r,x)+\frac\kappa
2\pa_2^2\til\lambda(r,x)=0,\label{distribution
covering}\end{equation} where $\pa_1$ and $\pa_2$ are partial
derivatives w.r.t.\ the first and second variable, respectively.

Let $\til\Lambda(p,x)=\int_0^x\til\lambda(p,s)ds$ for $p>0$ and
$x\in\R$. Then  for any $p>0$, $\til\Lambda(p,\cdot)$ is an odd and
increasing function, $\lim_{x\to\pm\infty}\til\Lambda(p,x)=\pm\frac
12$, and $\til\lambda(p,x)=\pa_2\til\Lambda(p,x)$. Thus for any
$r>0$ and $x\in\R$,
\begin{equation*}\pa_2(-\pa_1\til\Lambda(r,x)+\ha\HA_r(x)\pa_2\til\Lambda(r,x)
+\frac\kappa 2\pa_2^2\til\Lambda(r,x))=0.
\end{equation*}
Since $\til\Lambda(r,\cdot)$ is an odd function and $\ha\HA_r(0)=0$,
so
$$-\pa_1\til\Lambda(r,0)+\ha\HA_r(0)\pa_2\til\Lambda(r,0)
+\frac\kappa 2\pa_2^2\til\Lambda(r,0)=0.$$ Thus for any $r>0$ and
$x\in\R$, we have
\begin{equation}-\pa_1\til\Lambda(r,x)+\ha\HA_r(x)\pa_2\til\Lambda(r,x)
+\frac\kappa 2\pa_2^2\til\Lambda(r,x)=0.
\label{distribution-antiderivative}
\end{equation}
And we expect that for any $x\in\R\sem\{0\}$, $\lim_{r\to
0}\til\Lambda(r,x)\to\mbox{sign}\frac 12$. On the other hand, if
$\til\Lambda(r,x)$ satisfies (\ref{distribution-antiderivative}),
then $\til\lambda(r,x):=\pa_2\til\Lambda(r,x)$ satisfies
(\ref{distribution covering}).

 Let
$\lambda(r,x)=\sum_{k\in\Z}\til\lambda(r,x+2k\pi)$. Then
$\lambda(r,\cdot)$ has a period $2\pi$, and is the density function
of the distribution of the argument of $\lim_{t\to r}\beta(t)$,
where $\beta$ is a standard annulus SLE$_\kappa$ trace of modulus
$r$. So it satisfies $\int_{-\pi}^{\pi}\lambda(r,x)dx=1$. And
$\lambda(r,\cdot)$ is an even function for any $r>0$. Since
$\ha\HA_r$ has a period $2\pi$, so $\lambda(r,x)$ also satisfies
equation (\ref{distribution covering}). Let
$\Lambda(r,x)=\int_0^x\lambda(r,s)ds$. Then $\Lambda(r,x)$ satisfies
(\ref{distribution-antiderivative}). But $\Lambda(r,x)$ does not
satisfies $\lim_{x\to\pm\infty}\Lambda(r,x)=\pm 1$. Instead, we have
$\Lambda(r,x+2\pi)=\Lambda(r,x)+1$. In the case that $\kappa=2$, we
have some nontrivial solutions to
(\ref{distribution-antiderivative}). From Lemma 3.1 in \cite{Zhan},
we see $-\pa_r\HA_r+\HA_r\HA_r'+\HA_r''=0$, where the function
$\til\SA_r$ in \cite{Zhan} is the function $\HA_r$ here. From the
definition of $\ha\HA_r$, we may compute that
$-\pa_r\ha\HA_r+\ha\HA_r\ha\HA_r'+\ha\HA_r''=0$. Thus
$\Lambda_1(r,x)=\ha\HA_r(x)$ and $\Lambda_2(r,x)=r\HA_r(x)+x$
satisfy equation (\ref{distribution-antiderivative}). So
$\lambda_1(r,x)=\ha\HA_r'(x)$ and $\lambda_2(r,x)=r\HA_r'(x)+1$ are
solutions to (\ref{distribution covering}). In fact,
$\lambda_2(r,x)/(2\pi)$ is the distribution of the argument of the
end point of a Brownian Excursion in $\A_r$ started from $1$
conditioned to hit $\CC_r$. From Corollary 3.1 in \cite{Zhan}, this
is also the distribution of the argument of the limit point of a
standard annulus SLE$_2$ trace of modulus $r$. So we justified
equation (\ref{distribution covering}) in the case $\kappa=2$.

We may change variables in the following way. For $-\infty<s<0$, let
$\til{\mbox P}(s,y)=\til\Lambda(-\frac{\pi^2}s,-\frac\pi s y)$ and
${\mbox P}(s,y)=\Lambda(-\frac{\pi^2}s,-\frac\pi s y)$. Then for any
$s<0$, $\lim_{y\to\pm\infty}\til{\mbox P}(s,y)=\pm\frac 12$ and
${\mbox P}(s,y+2s)={\mbox P}(s,y)-1$. And we expect that
$\lim_{s\to-\infty}\til{\mbox P}(s,y)=\int_0^y\cosh(\frac
s2)^{-4/\kappa}ds/C(\kappa)$. Let $\GA_s(y)=i\HA_{-s}(iy-\pi)$ for
$s<0$ and $y\in\R$. From formula (\ref{change HA}), we may compute
that $\til\Lambda(r,x)$ ($\Lambda(r,x)$, resp.) satisfies equation
(\ref{distribution-antiderivative}) iff $\til{\mbox P}(s,y)$
(${\mbox P}(s,y)$, resp.) satisfies
\begin{equation}
-\pa_1\til{\mbox P}(s,y)+\GA_s(y)\pa_2\til{\mbox P}(s,y)+\frac\kappa
2\pa_2^2\til{\mbox P}(s,y)=0.\label{antiderivative after changing
variables}
\end{equation} From the
equation for $\HA_r$ and the definition of $\GA_s$, we have
$-\pa_s\GA_s+\GA_s\GA_s'+\GA_s''=0$. Thus ${\mbox
P}_1(s,y)=\GA_s(y)$ and ${\mbox P}_2(s,y)=s\GA_s(y)+y$ are solutions
to (\ref{antiderivative after changing variables}). In fact, ${\mbox
P}_1(s,y)$ corresponds to $-\Lambda_2(r,x)/\pi$, and ${\mbox
P}_2(s,y)$ corresponds to $-\pi\Lambda_1(r,x)$.

\section{Local Martingales for Annulus SLE$_4$ and SLE$_8$}
\subsection{Annulus SLE$_4$}
 Fix $\kappa=4$. Let $K_t$ and $\vphi_t$, $0\le t<p$, be
the annulus LE hulls and maps of modulus $p$, respectively, driven
by $\xi(t)=\sqrt\kappa B(t)$. Let $\beta(t)$, $0\le t<p$, be the
trace. For $r>0$, let $\TT^{(2)}_r(z)=\frac 12\SA_r(z^2)$ and
$\til\TT^{(2)}_r(z)=\frac 1i\TT^{(2)}_r(e^{iz})$. Solve the
differential equations:
$$\pa_t\psi_t(z)=\psi_t(z)\TT^{(2)}_{p-t}(\psi_t(z)/e^{i\xi(t)/2}),\,\,\,\psi_0(z)=z;$$
$$\pa_t\til\psi_t(z)=\til\TT^{(2)}_{p-t}(\til\psi_t(z)-\xi(t)/2),\,\,\,\til\psi_0(z)=z.$$
Let $P_2$ be the square map: $z\mapsto z^2$. Then we have
$P_2\circ\psi_t=\vphi_t\circ P_2$ and
$e^i\circ\til\psi_t=\psi_t\circ e^i$. Let $L_t:=P_2^{-1}(K_t)$ and
$\til L_t=(e^i)^{-1}(L_t)$. Then $\psi_t$ maps $\A_{p/2}\sem L_t$
conformally onto $\A_{(p-t)/2}$, and $\til\psi_t$ maps
$\St_{p/2}\sem\til L_t$ conformally onto $\St_{(p-t)/2}$. Since
$K_t=\beta(0,t]$, and $\beta$ is a simple curve in $\A_p$ with
$\beta(0)=1$, so $L_t$ is the union of two disjoint simple curves
opposite to each other, started from $1$ and $-1$, respectively. Let
$\alpha_\pm(t)$ denote the curve started from $\pm1$. Then
$\psi_{t}(\alpha_\pm(t))=e^i(\pm \xi(t)/2)$.

For each $r>0$, suppose $J_r$ is the conformal map from $\A_{r/2}$
onto $\{z\in\C:|\Imm z|<1\}\sem[-a_r,a_r]$ for some $a_r>0$ such
that $\pm 1$ is mapped to $\pm\infty$. This $J_r$ is symmetric
w.r.t.\ both $x$-axis and $y$-axis, i.e.,
$J_r(\lin{z})=\lin{J_r(z)}$, and $J_r(-z)=-J_r(z)$. And $\Imm J_r$
is the unique bounded harmonic function in $\A_{r/2}$ that satisfies
(i) $\Imm J_r\equiv \pm 1$ on the open arc of $\CC_0$ from $\pm1$ to
$\mp1$ in the ccw direction; and (ii) $\Imm J_r\equiv 0$ on
$\CC_{r/2}$. Let $\til J_r=J_r\circ e^i$.

\begin{Lemma} $-\pa_r\til J_r+\til J_r'\til\TT^{(2)}_r+\frac 12 \til
J_r''\equiv 0$ in $\til\A_{r/2}$.\end{Lemma} {\bf Proof.} Since
$\Imm \til J_r\equiv 0$ on $\R_{r/2}$, by reflection principle,
$\til J_r$ can be extended analytically across $\R_{r/2}$. And we
have $\Imm \til J_r'=\pa_x\Imm\til J_r\equiv 0$ and $\Imm \til
J_r''=\pa_x^2\Imm\til J_r\equiv 0$ on $\R_{r/2}$. From the equality
$\Imm \til J_r(x+ir/2)=0$, we have $\pa_r\Imm\til J_r+\pa_y\Imm\til
J_r/2\equiv 0$ on $\R_{r/2}$. On $\R_{r/2}$, note that $\Imm
\til\TT^{(2)}_r\equiv -1/2$, so
$$\Imm(\til J_r'\til\TT^{(2)}_r)=\Ree\til J_r'\Imm\til\TT^{(2)}_r+\Imm\til J_r'
\Ree\til\TT^{(2)}_r$$ $$=-1/2\Ree\til J_r'=-1/2\pa_y\Imm\til J_r
=\pa_r\Imm\til J_r.$$ Let $F_r:=-\pa_r\til J_r+\til
J_r'\til\TT^{(2)}_r+\frac 12 \til J_r''$. Then $\Imm F_r\equiv 0$ on
$\R_{r/2}$.

For any $k\in\Z$, we see that $\til J_r(z)$ is equal to
$(-1)^{k+1}\frac 2\pi \ln(z-k\pi)$ plus some analytic function for
$z\in\til\A_{r/2}$ near $k\pi$. So we may extend $\Ree\til J_r(z)$
harmonically across $\R\sem\{k\pi:k\in\Z\}$. Since $\Imm\til J_r$
takes constant value $(-1)^k$ on each interval $(k\pi,(k+1)\pi)$,
$k\in\Z$, we have $\Ree\til J_r(\lin z)=\Ree\til J_r(z)$. Moreover,
the following properties hold: $\pa_r \til J_r$ is analytic in a
neighborhood of $\R$, $\til J_r'$ and $\til J_r''$ are analytic in a
neighborhood of $\R\sem\{k\pi:k\in\Z\}$.

The fact that $\Imm\til J_r$ takes constant value $(-1)^k$ on each
$(k\pi,(k+1)\pi)$, $k\in\Z$, implies that $\Imm\pa_r\til J_r$,
$\Imm\til J_r'$ and $\Imm\til J_r''$ vanishes on
$\R\sem\{k\pi:k\in\Z\}$. Since $\Imm\til\TT^{(2)}_r$ also vanishes
on $\R\sem\{k\pi:k\in\Z\}$, so we compute $\Imm F_r\equiv 0$ on
$\R\sem\{k\pi:k\in\Z\}$.

From $J_r(\lin z)=\lin{J_r(z)}$, we find that $\til J_r(-\lin
z)=\lin{\til J_r(z)}$. So $\Ree \til J_r(z)=\Ree\til J_r(-\lin
z)=\Ree\til J_r(-z)$. This means that $\Ree \til J_r$ is an even
function, so is $\pa_r \til J_r$ and $\til J_r''$. And $\til J_r'$
is an odd function. Note that $\til \TT^{(2)}_r$ is an odd function,
so $F_r$ is an even function. Since $\til \TT^{(2)}_r(z)$ is equal
to $1/(2z)$ plus
 some analytic function for $z$ near $0$, so the pole of $F_r$ at $0$
 has order at most $2$. However, the coefficient of $1/z^2$ is equal
 to $2/\pi*1/2-1/2*2/\pi=0$. And
$0$ is not a simple pole of $F_r$ because $F_r$ is even. So $0$ is a
removable pole of $F_r$. Similarly, $\pi$ is also a removable pole
of $F_r$. Since $F_r$ has period $2\pi$, so every $k\pi$, $k\in\Z$,
is a removable pole of $F_r$. So $F_r$ can be extended analytically
across $\R$, and $\Imm F_r\equiv 0$ on $\R$. Thus  $\Imm F_r\equiv
0$ in $\St_{r/2}$, which implies that $F_r\equiv C$ for some
constant $C\in\R$.

Finally, $J_r(-z)=-J_r(z)$ implies that $\til J_r(z+\pi)=-\til
J_r(z)$. Since $\pi$ is a period of $\til \TT^{(2)}_r$, we compute
$F_r(z+\pi)=-F_r(z)$. So $C$ has to be $0$. $\Box$

\begin{Proposition} For any $z\in\A_{p/2}$,
$J_{p-t}(\psi_t(z)/e^{i\xi(t)/2})$ is a local martingale, from which
follows that $\Imm J_{p-t}(\psi_t(z)/e^{i\xi(t)/2})$ is a bounded
martingale.\label{right after}
\end{Proposition}
{\bf Proof.} Fix $z_0\in\St_{p/2}$, let
$Z_t:=\til\psi_t(z_0)-\xi(t)/2$, then
$$dZ_t=\til\TT^{(2)}_{p-t}(Z_t)dt-d\xi(t)/2.$$
Note that $\xi(t)/2=B(t)$. From Ito's formula and the last lemma, we
have
$$d \til J_{p-t}(Z_t)=-\pa_r\til J_{p-t}(Z_t)dt+\til
J_{p-t}'(Z_t)dZ_t + \frac 12 \til J_{p-t}''(Z_t)dt$$
$$=(-\pa_r\til J_{p-t}(Z_t)+\til J_{p-t}(Z_t)\til
\TT^{(2)}_{p-t}(Z_t)+\frac 12\til J_{p-t}''(Z_t))dt-\til
J_{p-t}'(Z_t)d\xi(t)/2=-\til J_{p-t}'(Z_t)d\xi(t)/2.$$ Thus $\til
J_{p-t}(Z_t)$, $0\le t<p$, is a local martingale. For any
$z\in\A_{p/2}$, there is $z_0\in\St_{p/2}$ such that $z=e^i(z_0)$.
Then
$$J_{p-t}(\psi_t(z)/e^{i\xi(t)/2})=J_{p-t}(\psi_t(e^i(z_0))/e^{i\xi(t)/2})
$$$$=J_{p-t}(e^i(\til\psi_t(z_0)-{\xi(t)/2}))=\til
J_{p-t}(\til\psi_t(z_0)-\xi(t)/2.$$ So
$J_{p-t}(\psi_t(z)/e^{i\xi(t)/2})$, $0\le t<p$, is a local
martingale. Since $|\Imm J_r(z)|\le 1$ for any $r>0$ and $z\in
\A_{r/2}$, so $\Imm J_{p-t}(\psi_t(z)/e^{i\xi(t)/2})$, $0\le t<p$,
is a bounded martingale. $\Box$

\vskip 3mm

 Let $h_t(z)=J_{p-t}(\psi_t(z)/e^{i\xi(t)/2})$. Then $h_t$
maps $\A_{(p-t)/2}\sem L_t$ conformally onto $\{z\in\C:|\Imm
z|<1\}\sem [-a_{p-t},a_{p-t}]$ so that $\alpha_\pm(t)$ is mapped to
$\pm\infty$. So $\Imm h_t$ is the unique bounded harmonic function
in $\A_{p/2}\sem L_t$ that vanishes on $\CC_{p/2}$, equals to $1$ on
the arc of $\CC_0$ from $1$ to $-1$ in the ccw direction and the
north side of $\alpha_+(0,t)$ and $\alpha_-(0,t)$, and equals to
$-1$ on the arc of $\CC_0$ from $-1$ to $1$ in the ccw direction and
the south side of $\alpha_+(0,t)$ and $\alpha_-(0,t)$.

We have another choice of $J_r$. Let $J_r$ be the conformal map of
$\A_{r/2}$ onto the strip $\{z\in\C:|\Imm z|<1\}\sem[-ib_r,ib_r]$
for some $b_r>0$ so that $\pm 1$ is mapped to $\pm\infty$. Then
$\Imm J_r$ is the bounded harmonic function in $\A_{r/2}$ determined
by the following properties: (i) $\Imm J_r\equiv \pm1$ on the open
arc of $\CC_0$ from $\pm1$ to $\mp1$ in the ccw direction; and (ii)
the normal derivative of $\Imm J_r$ vanishes on $\CC_{r/2}$. Let
$\til J_r=J_r\circ e^i$. The proposition and lemma in this
subsection still hold for $J_r$ and $\til J_r$ defined here. The
proofs are almost the same. The only difference is at the step when
we prove $\Imm F_r\equiv 0$ on $\R_{r/2}$. Here we have $\Ree\til
J_r'= \pa_y\Imm\til J_r$ vanishes on $\R_{r/2}$. Use an argument
similar to the proof of the lemma, we can show that $\Ree F_r'$
vanishes on $\R_{r/2}$. So $\pa_y\Imm F_r$ vanishes on $\R_{r/2}$.
Since $\Imm F_r$ vanishes on $\R$, $\Imm F_r$ has to vanish in
$\St_{r/2}$. Use $J_r(-z)=-J_r(z)$, we then conclude that $F_r\equiv
0$. If we let $h_t(z)=J_{p-t}(\psi_t(z)/e^{i\xi(t)/2})$, then $\Imm
h_t$ is the unique bounded harmonic function in $\A_{p/2}\sem L_t$
that satisfies the following properties: equals to $1$ on the arc of
$\CC_0$ from $1$ to $-1$ in the ccw direction, and the north side of
$\alpha_+(0,t)$ and $\alpha_-(0,t)$, equals to $-1$ on the arc of
$\CC_0$ from $-1$ to $1$ in the ccw direction, and the south side of
$\alpha_+(0,t)$ and $\alpha_-(0,t)$, and the normal derivatives
vanish on $\CC_{p/2}$.

Suppose $E$ is a doubly connected domain such that $0$ lies in the
bounded component of $\C\sem E$. Fix $v\in\pa_o E$, the outside
boundary component of $E$. Let $\beta(t)$, $0\le t<p$, be an
SLE$_4(E;v\to \pa_i E)$ trace, where $\pa_i E$ is the inside
boundary component of $E$. So $\beta(t)$ is  the image of a standard
annulus SLE$_4$ trace of modulus $p$ under the conformal map from
$(\A_p,1)$ onto $(E,v)$, where $p$ is the modulus of $E$. Let
$D=P_2^{-1}(E)$ and $\{v_+,v_-\}=P_2^{-1}(v)$.
$P_2^{-1}(\{\beta(t):0\le t<p\})$ is the union of two disjoint
simple curve started from $v_+$ and $v_-$, respectively. Let
$\alpha_\pm$ denote the curve started from $v_\pm$. Then $D$ is a
symmetric ($-D=D$) doubly connected domain, and
$\alpha_-(t)=-\alpha_+(t)$ for $0\le t<p$. Let
$D_t=D\sem\alpha_-([0,t])\sem\alpha_+([0,t])$. Let $\gamma_t^\pm$
denote the boundary arc of $\pa_o D_t$ from $\alpha_\pm(t)$ to
$\alpha_\mp(t)$ in the ccw direction. Then $\gamma_t^\pm$ contains a
boundary arc of $\pa_o D$, one side of $\alpha_+([0,t])$ and one
side of $\alpha_-([0,t])$. Let $H_t$ be the bounded harmonic
function in $D_t$ which has continuations at $\pa_i D$ and
$\gamma_t^\pm$ such that $H_t\equiv 0$ on $\pa_i D$ and
$H_t\equiv\pm 1$ on $\gamma_t^\pm$. By the definition of
SLE$_4(E;v\to \pa_i E)$ and conformal invariance of harmonic
functions, for any fixed $z_0\in D$, $H_t(z_0)$, $0\le t<p$, is a
bounded martingale. This $H_t$ corresponds to $\Imm h_t$ defined
right after Proposition \ref{right after}. We may replace the
condition $H_t\equiv 0$ on $\pa_i D$  by $\pa_{\bf n} H_t\equiv 0$
on $\pa_i D$. Then this $H_t$ corresponds to the $\Imm h_t$ defined
in the last paragraph. So for any fixed $z_0\in D$, it is still true
that  $H_t(z_0)$, $0\le t<p$, is a bounded martingale.

\subsection{Harmonic Explorers for Annulus SLE$_4$}
Let $D$ be a symmetric ($-D=D$) doubly connected subset of hexagonal
faces in the planar honeycomb lattice. Two faces of $D$ are
considered adjacent if they share an edge. Let $\pa_o D$ and $\pa_i
D$ denote the outside and inside component of $\pa D$, respectively.
Suppose $v_+$ and $v_-$ are vertices that lie on $\pa_o D$, and are
opposite to each other, i.e., $v_-=-v_+$. Suppose $\Ree v_+>0$. Then
$v_+$ and $v_-$ partition the boundary faces of $D$ near $\pa_o D$
into an "upper" boundary component, colored black, and a "lower"
boundary component, colored white. All other hexagons in $D$ are
uncolored.

Now we construct two curves $\alpha_+$ and $\alpha_-$ as follows.
Let $\alpha_\pm(0)=v_\pm$. Let $\alpha_\pm(1)$ be a neighbor vertex
of $\alpha_\pm(0)$ such that $[\alpha_\pm(0),\alpha_\pm(1)]$ is
shared by a white hexagon and a black hexagon. At time $n\in\N$, if
$\alpha_\pm(n)\not\in\pa_iD$, then $\alpha_\pm(n)$ is a vertex
shared by a black hexagon, a white hexagon, and an uncolored
hexagon, denoted by $f^n_\pm$. Let $H_n$ be the function defined on
faces, which takes value $1$ on the black faces, $-1$ on the white
faces, $0$ on faces that touch $\pa_i D$, and is discrete harmonic
at other faces of $D$. Then $H_n(f^n_-)=-H_n(f^n_+)$. We then color
$f_\pm$ black with probability equal to $(1+H_n(f^n_\pm))/2$ and
white with probability equal to $(1-H_n(f^n_\pm))/2$ such that
$f^n_+$ and $f^n_-$ are colored differently. Let $\alpha_\pm(n+1)$
be the unique neighbor vertex of $\alpha_\pm(n)$ such that
$[\alpha_\pm(n),\alpha_\pm(n+1)]$ is shared by a white hexagon and a
black hexagon. Increase $n$ by $1$, and iterate the above process
until $\alpha_+$ and $\alpha_-$ hit $\pa_i D$ at the same time. We
always have $\alpha_-(n)=-\alpha_+(n)$, $f^n_-=-f^n_+$, and
$H_n(-g)=-H_n(g)$.

From the construction, conditioned on $\alpha_\pm(k)$,
$k=0,1,\dots,n$, the expected value of $H_{n+1}(f^n_\pm)$ is equal
to $(1+H_n(f^n_\pm))/2-(1-H_n(f^n_\pm))/2=H_n(f^n_\pm)$. And if a
face $f$ is colored before time $n$, then its color will not be
changed after time $n$, so $H_{n+1}(f)=H_n(f)$. Since $H_{n+1}$ and
$H_n$ both vanish on the faces near $\pa_i D$, and are discrete
harmonic at all other uncolored faces at time $n+1$ and $n$, resp.,
so for any face $f$ of $D$, the conditional value of $H_{n+1}(f)$
w.r.t.\ $\alpha_\pm(k)$, $k=0,1,\dots,n$ is equal to $H_n(f)$. Thus
for any fixed face $f_0$ of $D$, $H_n(f_0)$ is a martingale.

If $n-1<t<n$, and $\alpha_\pm(n-1)$ and $\alpha_\pm(n)$ are defined,
let $\alpha_\pm(t)=(n-t)\alpha_\pm(n-1)+(t-(n-1))\alpha_\pm(n)$.
Then $\alpha_\pm$ becomes a curve in $D$. Let
$D_t=D\sem\alpha_+([0,t])\sem\alpha_-([0,t])$. Note that if the side
length of the hexagons is very small compared with the size of $D$,
then for any face $f$ of $D$, $H_n(f)$ is close to the value of
$\til H_n$ at the center of $f$, where $\til H_n$ is the bounded
harmonic function defined on $D_n$, which has a continuation to $\pa
D\sem\{v_+,v_-\}$ and the two sides of $\alpha_\pm([0,t))$ such that
$\til H_n\equiv 0$ on $\pa_i D$, and $\til H_n\equiv \pm 1$ on the
curve on $\pa_o D_n$ from $\alpha_\pm(n)$ to $\alpha_\mp(n)$ in the
ccw direction. From the last section, we may guess that the
distribution of $\alpha_\pm$ tends to that of the square root of an
annulus SLE$_4(P_2(D);P_2(v_\pm)\to\pa_i P_2(D))$ trace when the
mesh tends to $0$. If at each step of the construction of
$\alpha_\pm$, we let $ H_n$ be the function which is is equal to $1$
on the black faces, $-1$ on the white faces, and is discrete
harmonic at all other faces of $D$ including the faces that touch
$\pa_i D$, then we get a different pair of curves $\alpha_\pm$. If
the mesh is very small compared with the size of $D$, then for any
face $f$ of $D$, $H_n(f)$ is close to the value of $\til H_n$ at the
center of $f$, where $\til H_n$ is the bounded harmonic function
defined on $D_n$, which has a continuation to $\pa D\sem\{v_+,v_-\}$
and the two sides of $\alpha_\pm([0,t))$ such that $\pa_{\bf n}\til
H_n\equiv 0$ on $\pa_i D$, and $\til H_n\equiv \pm 1$ on the curve
on $\pa_o D_n$ from $\alpha_\pm(n)$ to $\alpha_\mp(n)$ in the ccw
direction. So we also expect the law of $\alpha_\pm$ constructed in
this way tends to that of the square root of an annulus
SLE$_4(P_2(D);P_2(v_\pm)\to\pa_i P_2(D))$ trace when the mesh tends
to $0$.

\subsection{Annulus SLE$_8$}

Fix $\kappa=8$. Let $K_t$ and $\vphi_t$, $0\le t<p$, be the annulus
LE hulls and maps, respectively, of modulus $p$, driven by
$\xi(t)=\sqrt\kappa B(t)$. For $r>0$, let $\TT^{(4)}_r(z)=\frac
14\SA_r(z^4)$ and $\til\TT^{(4)}_r(z)=\frac 1i\TT^{(4)}_r(e^{iz})$.
Solve the differential equations:
$$\pa_t\psi_t(z)=\psi_t(z)\TT^{(4)}_{p-t}(\psi_t(z)/e^{i\xi(t)/4}),\,\,\,\psi_0(z)=z;$$
$$\pa_t\til\psi_t(z)=\til\TT^{(4)}_{p-t}(\til\psi_t(z)-\xi(t)/4),\,\,\,\til\psi_0(z)=z.$$
Let $P_4$ be the map: $z\mapsto z^4$. Then we have
$P_4\circ\psi_t=\vphi_t\circ P_4$ and
$e^i\circ\til\psi_t=\psi_t\circ e^i$. Let $L_t:=P_4^{-1}(K_t)$ and
$\til L_t=(e^i)^{-1}(L_t)$. Then $\psi_t$ maps $\A_{p/4}\sem L_t$
conformally onto $\A_{(p-t)/4}$, and $\til\psi_t$ maps
$\St_{p/4}\sem\til L_t$ conformally onto $\St_{(p-t)/4}$. Let $G_r$
map $\A_{r/4}$ conformally onto $\{z\in\C:|\Ree z|+|\Imm
z|<1\}\sem[-a_r,a_r]$ for some $a_r>0$ such that $\pm 1$ and $\pm i$
are fixed.

\begin{Proposition} For any $z\in\A_{r/4}$,
$G_{p-t}(\psi_t(z)/e^{i\xi_t/4})$ is a bounded martingale.
\end{Proposition}
{\bf Proof.} Let $\til G_r:=G_r\circ e^i$. For any $z\in\A_{p/4}$,
there is $w\in\St_{p/4}$ such that $z=e^i(w)$. Then
$$G_{p-t}(\psi_t(z)/e^{i\xi(t)/4})=G_{p-t}(\psi_t(e^{iw})/e^{i\xi(t)/4})=\til
G_{p-t}(\til\psi_t(w)-\xi(t)/4).$$ To prove this proposition, it
suffices to show that for any $w\in\St_{p/4}$, $\til
G_{p-t}(\til\psi_t(w)-\xi(t)/4)$ is a local martingale. Let
$Z_t=\til\psi_t(w)-\xi(t)/4$, then
$$dZ_t=\til\TT^{(4)}_{p-t}(Z_t)dt-dB(t)/\sqrt 2.$$
Thus by Ito's formula,
$$d\til G_{p-t}(\psi_t(w)-\xi(t)/4)=-\pa_r\til G_{p-t}(Z_t)dt+\til G_{p-t}'(Z_t)dZ_t
+\frac 12\til G_{p-t}''(Z_t)\frac{dt}2$$
$$=(-\pa_r\til G_{p-t}(Z_t)+\til G_{p-t}'(Z_t)\til\TT^{(4)}_{p-t}(Z_t)+\frac
14\til G_{p-t}''(Z_t))dt-\til G_{p-t}'(Z_t)dB(t)/\sqrt 2.$$ So it
suffices to prove the following lemma.

\begin{Lemma} $-\pa_r\til G_r+\til G_r'\til \TT^{(4)}_r+\frac 14 \til G_r''\equiv 0$
in $\St_{r/4}$. \end{Lemma} {\bf Proof.} Let $F_r$ be the left-hand
side. Let $Q_r(z):=i(\til G_r(z)-1)^2$. Note that $\til G_r$ maps
$[0,\pi/2]$ and $[-\pi/2,0]$ onto the line segments $[1,i]$ and
$[-i,1]$, respectively. Thus $Q_r(z)\to\R$ as $z\in\St_{r/4}$ and
$z\to (-\pi/2,\pi/2)$. By reflection principle, $Q_r$ can be
extended to an analytic function in a neighborhood of
$(-\pi/2,\pi/2)$, and $Q_r(\lin{z})=\lin{Q_r(z)}$. Since
$G_r(\lin{z})=\lin{G_r(z)}$, so $\til G_r(-\lin z)=\lin{\til
G_r(z)}$. It follows that $Q_r(-\lin z)=-\lin{Q_r(z)}$. So we have
$Q_r(-z)=-Q_r(z)$, and the Taylor expansion of $Q_r$ at $0$ is
$\sum_{n=0}^\infty a_{n,r}z^{2n+1}$. Thus $\til
G_r(z)=1+\sum_{n=0}^\infty c_{n,r}z^{2n+1/2}$ for $z$ near $0$. So
$\pa_r\til G_r(z)=O(z^{1/2})$ for $z$ near $0$, $\til
G_r'(z)=1/2c_{1,r}z^{-1/2}+O(z^{3/2})$, and $\til
G_r''(z)=-1/4c_{1,r}z^{-3/2}+O(z^{1/2})$. Since
$\til\TT^{(4)}_r(z)=1/(8z)+O(z)$ near $0$, so $\til G_r'(z)\til
\TT^{(4)}_r(z)=1/16 c_{1,r}z^{-3/2}+O(z^{1/2})$. Then we compute
$F_r(z)=O(z^{1/2})$ near $0$. Similarly,
$F_r(z)=O((z-k\pi/2)^{1/2})$ for $z$ near $k\pi/2$, $k\in\Z$.

For $z\in(k\pi,(k+1/2)\pi)$, $k\in\Z$, $\til G_r(z)\in
(-1)^k+(1-i)\R$. So $F_r(z)\in(1-i)\R$ for $z\in(k\pi,(k+1/2)\pi)$,
$k\in\Z$. Similarly, $F_r(z)\in (1+i)\R$ for
$z\in((k-1/2)\pi,k\pi)$, $k\in\Z$. Since $\til G_r$ takes real
values on $\R_{r/4}$, so $F_r$ also takes real values on $\R_{r/4}$.
 Let $V_r=\Imm F_r$, then $V_r\equiv 0$ on $\R_{r/4}$, and for $k\in\Z$,
 $\pa_x V+\pa_y V\equiv 0$ on $(k\pi,(k+1/2)\pi)$ and $\pa_x V-\pa_y V\equiv 0$
on $((k-1/2)\pi,k\pi)$, $k\in\Z$. And $V_r(z)\to 0$ as
$z\in\St_{r/4}$ and $z\to k\pi/2$, $k\in\Z$. Since $\til G_r$ and
$\til\TT^{(4)}_r$ have period $2\pi$, so does $F_r$. Thus $|V_r|$
attains its maximum in $\lin{\St_{r/4}}$ at some
$z_0\in\R\cup\R_{r/4}$. If $z_0\in\R_{r/4}$ or $z_0=k\pi/2$ for some
$k\in\Z$, then $V_r(z_0)=0$, and so $V_r$ vanishes in $\St_{r/4}$.
Otherwise, either $z_0\in (k\pi,(k+1/2)\pi)$ or
$z_0\in((k-1/2)\pi,k\pi)$ for some $k\in\Z$. In either cases, we
have $\pa_xV_r(z_0)=0$, so $\pa_y V_r(z_0)=0$ too. Thus
$F_r'(z_0)=0$. If $F_r$ is not constant in $\St_{r/4}$, then
$F_r(z_0)=1+a_m(z-z_0)^m+O((z-z_0)^{m+1}$ for $z$ near $z_0$. Then
it is impossible that $\Imm F_r(z_0)\ge \Imm F_r(z)$ for all
$z\in\{|z-z_0|<\eps,\Imm z\ge \Imm z_0\}$ or $\Imm F_r(z_0)\le \Imm
F_r(z)$ for all $z\in\{|z-z_0|<\eps,\Imm z\ge \Imm z_0\}$. This
contradiction shows that $F_r$ has to be constant in $\St_{r/4}$.
Since $F_r(z)\to 0$ as $z\to 0$, so this constant is $0$. We again
conclude that $V_r$ has to vanish in $\St_{r/4}$. $\Box$

\section{Annulus SLE$_{8/3}$ and the Restriction Property}
In this section, we fix $\kappa=8/3$ and $\alpha=5/8$. Let $\vphi_t$
and $K_t$, $0\le t<p$, be the annulus LE maps and hulls of modulus
$p$, driven by $\xi(t)=\sqrt{\kappa}B(t)$, $0\le t<p$. Let
$\til\vphi_t$ and $\til K_t$, $0\le t<p$, be the corresponding
annulus LE maps and hulls in the covering space. Let $A\ne\emptyset$
be a hull in $\A_p$ w.r.t. $\CC_p$ (i.e., $\A_p\sem A$ is a doubly
connected domain whose one boundary component is $\CC_p$) such that
$1\not\in\lin A$. So there is $t>0$ such that $K_t\cap A=\emptyset$.
Let $T_A$ be the biggest $T\in(0,p]$ such that
 for $t\in[0,T)$, $K_t\cap A=\emptyset$. Let $\vphi_A$ be the conformal map from $\A_p\sem A$
onto $\A_{p_0}$ such that $\vphi_A(1)=1$, where $p_0$ is equal to
the modulus of $\A_p\sem A$. Let $K'_t=\vphi_A(K_t)$, $0\le t<T_A$.
Let $h(t)$ equal $p_0$ minus the modulus of $\A_{p_0}\sem K'_t$.
Then $h$ is a continuous increasing function with $h(0)=0$. So $h$
maps $[0,T_A)$ onto $[0,S_A)$ for some $S_A\in(0,p_0]$. From
Proposition 2.1 in \cite{Zhan}, $L_s=K_{h^{-1}(s)}$, $0\le s<S_A$,
are the annulus LE hulls of modulus $p_0$, driven by some real
continuous function, say $\eta(s)$. Let $\psi_s$, $0\le s<S_A$, be
the corresponding annulus LE maps. Let $\til\psi_s$ and $\til L_s$,
$0\le s<S_A$, be the annulus LE maps and hulls, respectively, in the
covering space.

Let $f_t=\psi_{h(t)}\circ\vphi_A\circ\vphi_t^{-1}$ and
$A_t=\vphi_t(A)$. Then for $0\le t<T_A$,
$e^i(\xi(t))\not\in\lin{A_t}$, and $f_t$ maps $(\A_{p-t}\sem
A_t,\CC_{p-t})$ conformally onto $(\A_{p_0-h(t)},\CC_{p_0-h(t)})$.
And for any $z_0\in\CC_0\sem\lin{A_t}$, if $z\in\A_{p-t}\sem A_t$
and $z\to z_0$, then $f(z)\to\CC_0$. Thus $f_t$ can be extended
analytically across $\CC_0$ near $e^i(\xi(t))$. A proof similar to
those of Lemma 2.1 and 2.2 in \cite{Zhan} shows that
$f_t(e^i(\xi(t)))=e^i(\eta(h(t)))$, and
$h'(t)=|f_t'(e^i(\xi(t)))|^2$.

Let $\til\vphi_A$ be such that $e^i\circ\til\vphi_A=\vphi_A\circ
e^i$ and $\til\vphi_A(0)=0$. Let $\til
f_t=\til\psi_{h(t)}\circ\til\vphi_A \circ\til\vphi_t^{-1}$. Then
$e^i\circ \til f_t=f_t\circ e^i$, and so $e^i\circ\til
f_t(\xi(t))=e^i(\eta(h(t)))$. Thus $\til
f_t(\xi(t))=\eta(h(t))+2k\pi$ for some $k\in\Z$. Now we replace
$\eta(s)$ by $\eta(s)+2k\pi$. Then $\eta(s)$, $0\le s<S_A$, is still
a driving function of $L_s$, $0\le s<S_A$. And we have $\til
f_t(\xi(t))=\eta(h(t))$. Moreover, we have $h'(t)=\til
f_t'(\xi(t))^2$.

Let $\til A=(e^i)^{-1}(A)$ and $\til A_t=(e^i)^{-1}(A_t)$.  For any
$t\in[0,T_A)$, and $z\in\St_p\sem \til A\sem \til K_t$, we have
$\til f_t\circ\til\vphi_t(z)=\til\psi_{h(t)}\circ \til \vphi_A(z)$.
Taking the derivative w.r.t.\ $t$, we compute
$$\partial_t \til f_t(\til\vphi_t(z))+\til f_t'(\til\vphi_t(z))
\HA_{p-t}(\til\vphi_t(z)-\xi(t))=\til
f_t'(\xi(t))^2\HA_{p_0-h(t)}(\til f_t(\vphi_t(z))-\til
f_t(\xi(t))).$$

Since $\til A_t=\til\vphi_t(A)$ for $0\le t<T_A$, so for any
$t\in[0,T_A)$, and $w\in\St_{p-t}\sem \til A_t$, we have
$\til\vphi_t^{-1}(w)\in\St_p\sem \til A\sem\til K_t$. Thus
\begin{equation} \partial_t\til
f_t(w)=\til f_t'({\xi(t)})^2\HA_{p_0-h(t)}({\til f_t(w)-\til
f_t(\xi(t))})- \til f_t'(w)\HA_{p-t}({w-\xi(t)}).
\label{reseqn}\end{equation}

Recall that
$$\HA_r(z)=-i\lim_{M\to\infty}\sum_{k=-M}^M\frac{e^{2kr}+e^{iz}}
{e^{2kr}-e^{iz}}=\cot(z/2)+\sum_{k=1}^\infty-i(\frac{e^{2kr}+e^{iz}}
{e^{2kr}-e^{iz}}+\frac{e^{-2kr}+e^{iz}} {e^{-2kr}-e^{iz}})$$
$$=\cot(z/2)+\sum_{k=1}^\infty\frac{2\sin(z)}{\cosh(2kr)-\cos(z)}.$$
Let $$S_r=\sum_{k=1}^\infty \frac{2}{\cosh(2kr)-1}=\sum_{k=1}^\infty
\frac{1}{\cosh^2(kr)}.$$ Then the Laurent sires  expansion of
$\HA_r$ at $0$ is $\HA_r(z)=2/z+(S_r- 1/6)z+O(z^2)$.

Apply the following power series expansions:
$$\HA_{r}(z)=2/z+O(z);$$
$$\til f_t'(w)=\til f_t'(\xi(t))+\til f_t''(\xi(t))(w-\xi(t))+O((w-\xi(t))^2);$$
$$\til f_t(w)-\til f_t(\xi(t))=\til f_t'(\xi(t))(w-\xi(t))+\frac{\til f_t''(\xi(t))}2
(w-\xi(t))^2+O((w-\xi(t))^3).$$ After some straightforward
computation and letting $w\to\xi(t)$, we get

\begin{Lemma} $\pa_t \til f_t(\xi(t))=-3\til f_t''(\xi(t))$. \label{loc}
\end{Lemma}

Now differentiate equation (\ref{reseqn}) with respect to $w$. We
get
\begin{equation*} \partial_t\til
f_t'(w)=\til f_t'({\xi(t)})^2\til f_t'(w)\HA_{p_0-h(t)}'({\til
f_t(w)-\til f_t(\xi(t))})$$$$- \til
f_t''(w)\HA_{p-t}({w-\xi(t)})-\til f_t'(w)\HA_{p-t}'({w-\xi(t)}).
\end{equation*}

 Apply the previous power series
expansions and the following expansions:
$$\HA_r(z)=2/z+(S_r-1/6)z +O(z^2);$$
$$\HA_r'(z)=-{2}/{z^2}+(S_r-1/6)+O(z);$$
$$\til f_t''(w)=\til f_t''(\xi(t))+\til f_t'''(\xi(t))(w-\xi(t))+O((w-\xi(t))^2);$$
$$\til f_t'(w)=\til f_t'(\xi(t))+\til f_t''(\xi(t))(w-\xi(t))+\frac{\til f_t'''(\xi(t))}2
(w-\xi(t))^2+O((w-\xi(t))^3);$$
$$\til f_t(w)-\til f_t(\xi(t))=\til f_t'(\xi(t))(w-\xi(t))+\frac{\til f_t''(\xi(t))}2
(w-\xi(t))^2$$$$+\frac{\til
f_t'''(\xi(t))}6(w-\xi(t))^3+O((w-\xi(t))^4).$$ After some long but
straightforward computation and letting $w\to\xi(t)$, we get

\begin{Lemma} $$\frac{\pa_t\til f_t'(\xi(t))}{\til f_t'(\xi(t))}=\frac
12\left(\frac{\til f_t''(\xi(t))}{\til f_t'(\xi(t))}\right)^2-\frac
43\frac{\til f_t'''(\xi(t))}{\til f_t'(\xi(t))}$$$$ +\til
f_t'(\xi(t))^2(S_{p_0-h(t)}-1/6)-(S_{p-t}-1/6).$$

\end{Lemma}

From Ito's formula and the above lemma, we have
$$d\til f_t'(\xi(t))=\pa_t\til f_t'(\xi(t))dt+\til f_t''(\xi(t))d\xi(t)+\frac\kappa
2\til f_t'''(\xi(t))dt$$
$$=\til f_t''(\xi(t))d\xi(t)+\til f_t'(\xi(t))\left(\frac
12\left(\frac{\til f_t''(\xi(t))}{\til
f_t'(\xi(t))}\right)^2+\left(\frac\kappa 2-\frac 43\right)\frac{\til
f_t'''(\xi(t))}{\til f_t'(\xi(t))}\right.$$$$ +\left.\til
f_t'(\xi(t))^2(S_{p_0-h(t)}-1/6) -(S_{p-t}- 1/6)\right)dt.$$

Thus
$$d\til f_t'(\xi(t))^\alpha=\alpha \til f_t'(\xi(t))^{\alpha-1}d\til f_t'(\xi(t))
+\alpha(\alpha-1)\til f_t'(\xi(t))^{\alpha-2}\frac{\kappa}2\til
f_t''(\xi(t))^2dt$$
$$=\alpha \til f_t'(\xi(t))^\alpha\left(\frac{d\til f_t'(\xi(t))}{\til f_t'(\xi(t))}
+(\alpha-1)\frac{\kappa}2\left(\frac{\til f_t''(\xi(t))}{\til
f_t'(\xi(t))}\right)^2dt\right)$$
$$=\alpha \til f_t'(\xi(t))^\alpha\left(\frac{\til f_t''(\xi(t))}{\til f_t'(\xi(t))}d\xi(t)+\left(\left(\frac
12+(\alpha-1)\frac{\kappa}2\right)\left(\frac{\til
f_t''(\xi(t))}{\til f_t'(\xi(t))}\right)^2\right.\right.
$$$$+\left.\left.\left(\frac\kappa 2-\frac
43\right)\frac{\til f_t'''(\xi(t))}{\til f_t'(\xi(t))}+\til
f_t'(\xi(t))^2\left(S_{p_0-h(t)}-\frac16\right) -\left(S_{p-t}-
\frac16\right)\right)dt\right)$$
$$=\alpha \til f_t'(\xi(t))^\alpha\left(
\frac{\til f_t''(\xi(t))}{\til
f_t'(\xi(t))}d\xi(t)+\left(h'(t)\left(S_{p_0-h(t)}-
\frac16\right)-\left(S_{p-t}- \frac16\right)\right)dt\right).$$ The
last equality uses $\kappa=8/3$, $\alpha=5/8$, and $h'(t)=\til
f_t'(\xi(t))^2$.

Now we have the following theorem.

\begin{Theorem} $$M_t=\til f_t'(\xi(t))^{5/8}\exp\left(-\frac{5}{8}
\int^{p-t}_{p_0-h(t)}\left(S_r-\frac16\right)dr\right),$$ $0\le
t<T_A$, is a bounded martingale.\label{bddmtgl}
\end{Theorem}
{\bf Proof.} From the above computation and Ito's formula, we see
tht $M_t$, $0\le t<T_A$, is a local martingale.

Since $f_t$ maps $\A_{p-t}\sem A_t$ conformally onto
$\A_{p_0-h(t)}$, so by the comparison principle of extremal length,
the modulus of $\A_{p_0-h(t)}$ is not bigger than that of
$\A_{p-t}$. Thus $p_0-h(t)\le p-t$. Since $S_r>0$ for any $r>0$, so
\begin{equation}\exp\left(-\frac 58\int_{p_0-h(t)}^{p-t}\left(S_r-\frac
16\right)dr\right)
\le\exp\left(\frac{5}{48}((p-t)-(p_0-h(t))\right)\le
\exp\left(\frac{5p}{48}\right).\label{exponential bounded}
\end{equation}

Let $g_t=f_t^{-1}$, $\til g_t=\til f_t^{-1}$. Then $g_t\circ
e^i=e^i\circ\til g_t$. And $g_t$ maps $\A_{p_0-h(t)}$ conformally
onto $\A_{p-t}\sem A_t$. Now $-\ln(g_t(z)/z)$ is a bounded analytic
function defined in $\A_{p_0-h(t)}$, $\Ree
(-\ln(g_t(z)/z))\to(p-t)-(p_0-h(t))$ as $z\to\CC_{p_0-h(t)}$, and
any subsequential limit of $\Ree (-\ln(g_t(z)/z))$ as $z\to\CC_0$ is
nonnegative. Thus there are some $C\in\R$ and a positive measure
$\mu_t$ supported by $\CC_0$ of total mass $(p-t)-(p_0-h(t))$ such
that for any $z\in\A_{p_0-h(t)}$,
\begin{equation}-\ln(g_t(z)/z)=\int_{\CC_0}\SA_{p_0-h(t)}(z/\theta)d\mu_t(\theta)+iC.
\label{Annulus integral}\end{equation}
 For
any $w\in\St_{p_0-h(t)}$, we have $e^i(w)\in\A_{p_0-h(t)}$,
$\ln(g_t(e^i(w)))=i\til g_t(w)$, and $\ln(z)=iw$, so
$$-i(\til g_t(w)-w)=\int_{\CC_0}\SA_{p_0-h(t)}(e^i(w)/\theta)d\mu_t(\theta)+iC.$$
If $\til\mu_t$ is a measure on $\R$ that satisfies
$\mu_t=\til\mu_t\circ (e^i)^{-1}$, then  for any
$w\in\St_{p_0-h(t)}$,
$$\til g_t(w)-w=\int_\R i\SA_{p_0-h(t)}\circ
e^i(w-x)d\til\mu_t(x)-C=\int_\R-\HA_{p_0-h(t)}(w-x)d\til\mu_t(x)-C.$$
Taking derivative w.r.t.\ $w$, we have \begin{equation}\til
g_t'(w)-1=\int_\R-\HA_{p_0-h(t)}'(w-x)d\til\mu_t(x).
\label{difference integral}
\end{equation}

 From equation (\ref{change HA}) and the definition of $\HA_r$, we have
$$\HA_r(z)=i\frac\pi r\HA_{\pi^2/r}(i\frac\pi r z)-\frac{z}r
=\frac\pi r\lim_{M\to\infty}\sum_{k=-M}^M\frac{e^{2k\pi^2/r}+e^{-\pi
z/r}}{e^{2k\pi^2/r}-e^{-\pi z/r}}-\frac zr.$$ Thus
\begin{equation}\HA_r'(z)=\frac{\pi^2}{r^2}\sum_{k=-\infty}^\infty\frac{-2e^{\pi
z/r}e^{2k\pi^2/r}}{(e^{\pi z/r}-e^{2k\pi^2/r})^2}-\frac
1r.\label{HA-prime}\end{equation} So for $z\in\R$, we have
$\HA_r'(z)<0$. Apply this to equation (\ref{difference integral}).
We get $\til g_t'(\eta(h(t)))>1$. Thus $\til f_t'(\xi(t))\in(0,1)$.
Then from equation (\ref{exponential bounded}), we have
$$0\le M_t\le\exp\left(-\frac 58\int_{p_0-h(t)}^{p-t}
\left(S_r-\frac
16\right)dr\right)\le\exp\left(\frac{5p}{48}\right).$$ Since $M_t$,
$0\le t<T_A$, is uniformly bounded, so it is a bounded martingale.
$\Box$

\vskip 3mm

Now suppose that $A$ is a smooth hull, i.e., there is a smooth
simple closed curve $\gamma:[0,1]\to\A_p\cup\CC_0$ with
$\gamma((0,1))\subset\A_p$ and $\gamma(0)\ne\gamma(1)\in\CC_0$, and
$A$ is bounded by $\gamma$ and an arc on $\CC_0$ between $\gamma(0)$
and $\gamma(1)$.

If $T_A<p$, a proof similar to Lemma 6.3 in \cite{LSW-8/3} shows
that $\til f_t'(\xi(t))\to 0$ as $t\to T_A$. Thus $M_t\to 0$ as
$t\to T_A$ on the event that $T_A<p$. From now on, we suppose
$T_A=p$. Then $K_t$ approaches $\CC_p$ as $t\to p$ and is uniformly
bounded away from $A$. Then the modulus of $\A_p\sem K_t\sem A$
tends to $0$ as $t\to p$. Thus $p_0-h(t)\to 0$ as $t\to p$. So
$S_A=p_0$. Now $A_t=\vphi_t(A)$ is bounded by
$\gamma_t=\vphi_t(\gamma)$ and an arc on $\CC_0$ between
$\gamma_t(0)$ and $\gamma_t(1)$. So $A_t$ is also a smooth hull.
Thus $f_t$ and $g_t$ both extend continuously to the boundary of the
definition domain. And $f_t$ maps $\gamma_t$ to an arc on $\CC_0$.
Let $I_t$ denote this arc. Since $-\ln(g_t(z)/z)$ also extends
continuously to $\CC_0$, so the measure $\mu_t$ in equation
(\ref{Annulus integral}) satisfies
$$d\mu_t(z)=-\Ree \ln(g_t(z)/z)/(2\pi)d{\bf m}(z)=-\ln|g_t(z)|/(2\pi)d{\bf m},$$
where $\bf m$ is the Lebesgue measure on $\CC_0$ (of total mass
$2\pi$). Since $\ln|g_t(z)|=0$ for $z\in\CC_0\sem I_t$, so $\mu_t$
is supported by $I_t$. Let $\til\gamma$ be a continuous curve such
that $\gamma=e^i\circ\til\gamma$. Let
$\til\gamma_t=\til\vphi_t(\til\gamma)$ and $\til I_t=\til
f_t(\gamma_t)$. Then $e^i(\til\gamma_t)=\gamma_t$, $e^i(\til
I_t)=I_t$, and $\til I_t$ is a real interval. Let $\til\mu_t$ be a
measure supported by $\til I_t$ that satisfies
$d\til\mu_t(z)=\Imm\til g_t(z)/(2\pi)d{\bf m_\R}$ for $z\in\til
I_t$, where ${\bf m}_\R$ is the Lebesgue measure on $\R$. Since
$-\ln|g_t(e^i(z))|=\Imm\til g_t(z)$, so
$\mu_t=\til\mu_t\circ(e^i)^{-1}$. Thus equation (\ref{difference
integral}) holds for this $\til\mu_t$.

Now $\til\vphi_t$ maps $\St_p\sem \til K_t$ conformally onto
$\St_{p-t}$. Let $\Sigma_t$ be the union of $\St_p\sem \til K_t$,
its reflection w.r.t.\ $\R$, and $\R\sem\lin{\til K_t}$. By Schwarz
reflection principle, $\vphi_t$ extends analytically to $\Sigma_t$,
and maps $\Sigma_t$ conformally into $\{z\in\C:|\Imm z|<p-t\}$. For
every $z\in A$, the distance from $z$ to the boundary of $\Sigma_t$
is at least $d_0=\min\{p,dist(A,K_p)\}>0$, and the distance from
$\vphi_t(z)$ to the boundary of $\{z\in\C:|\Imm z|<p-t\}$ equals to
$p-t$. By Koebe's $1/4$ theorem, $|\vphi_t'(z)|\le 4(p-t)/d_0$. Let
$H=\max\{\Imm\til \gamma(u):u\in[0,1]\}$.  Since
$\til\gamma_t=\til\vphi_t\circ\til\gamma$, so
$H_t:=\max\{\Imm\til\gamma_t(u):u\in[0,1]\}\le{4(p-t)}H/d_0$.  A
proof similar as above shows that for any $z\in \til I_0$,
$|\til\psi_{h(t)}(z)|\le 4(p_0-h(t))/d_1$ for some $d_1>0$. Since
\begin{equation}\til I_t=\til
f_t(\til\gamma_t)=\til\psi_{h(t)}\circ\til\vphi_A(\til\gamma)
=\til\psi_{h(t)}(\til I_0),\label{It-I0}\end{equation} so $|\til
I_t|\le 4(p-t)/d_1|\til I_0|$. Thus $|\mu_t|=|\til\mu_t|\le
H_t|I_t|\le 16(p-t)(p_0-h(t))H|\til I_0|/(d_0d_1)$. Let
$C_0=16H|\til I_0|/(d_0d_1)$, then
\begin{equation}(p-t)-(p_0-h(t))=|\mu_t|=|\til\mu_t|\le
C_0(p-t)(p_0-h(t)).\label{estimation of mut}\end{equation}
 Thus \begin{equation}(p_0-h(t))/(p-t)\ge
1-C_0(p_0-h(t)).\label{p0-h(t)}\end{equation}

Since $\til\mu_t$ is supported by $\til I_t$, so from equation
(\ref{difference integral}) we have
$$\til g_t'(\eta(h(t)))-1=\int_{\til
I_t}-\HA'_{p_0-h(t)}(\eta(h(t))-x)d\til\mu_t(x).$$ Let
$\til\alpha(t)=\til\vphi_t^{-1}(\xi(t))$. Then $\til\alpha(t)$ is a
simple curve, and
$\alpha(t)=e^i(\til\alpha(t))=\vphi_t^{-1}(e^i(\xi(t)))$ is an
annulus SLE$_{8/3}$ trace. So $K_t=\alpha((0,t])$ for any $t\ge 0$.
Thus $\til K_t=\cup_{k\in\Z}(\til\alpha((0,t])+2k\pi)$. Let
$\til\beta(s)=\til\vphi_A(\til\alpha(h^{-1}(s)))$ for $0\le s<p_0$.
Since $\til L_{h(t)}=\til\vphi_A(K_t)$ and
$\til\vphi_A(z+2k\pi)=\til\vphi_A(z)+2k\pi$, so $\til
L_{h(t)}=\cup_{k\in\Z}(\til\beta((0,h(t)])+2k\pi)$. Now we compute
$$\til\psi_{h(t)}(\til\beta(h(t)))=\til\psi_{h(t)}\circ\til\vphi_A\circ\til\vphi_{t}^{-1}
(\xi(t))=\til f_{t}(\xi(t))=\eta(h(t)).$$ Thus $\til\psi_s$ maps the
left and right side of $\beta((0,h(t)))$ to intervals
$(b_-(t),\eta(h(t)))$ and $(\eta(h(t)),b_+(t))$, respectively, for
some $b_-(t)<\eta(h(t))<b_+(t)$. Therefore $\til\psi_{h(t)}$ maps
the $\til L_{h(t)}$ to $\cup_{k\in\Z}(b_-(t)+2k\pi,b_+(t)+2k\pi)$.
From equation (\ref{It-I0}), we have
$\cup_{k\in\Z}(l(t)+2k\pi,r(t)+2k\pi)\cap\til I_{t}=\emptyset$. So
for any $x\in\til I_{t}$ and $k\in\Z$,
$|x-(\eta(s)+2k\pi)|\ge\min\{\eta(h(t))-b_-(t),b_+(t)-\eta(h(t))\}$.

As $t\to p$,  $\til\beta(h(t))$ approaches to a point on $\R_{p_0}$,
so the extremal distance between the left side of
$\til\beta((0,h(t))$ and $\R_{p_0}$ in $\St_{p_0}\sem \til L_{h(t)}$
tends to $0$. Since $\til\vphi_t$ maps $(\St_{p_0}\sem \til
L_{h(t)},\R_{p_0})$ conformally onto
$(\St_{p_0-h(t)},\R_{p_0-h(t)})$, and the left side of
$\til\beta((0,h(t))$ is mapped to $(b_-(t),\eta(h(t)))$, so the
extremal distance between $(b_-(t),\eta(h(t)))$ and $\R_{p_0-h(t)}$
in $\St_{p_0-h(t)}$ tends to $0$ as $t\to p$ by the conformal
invariance property of extremal length. Thus
$(\eta(h(t))-b_-(t))/(p_0-h(t))\to+\infty$ as $t\to p$. Similarly,
$(b_+(t)-\eta(h(t)))/(p_0-h(t))\to+\infty$ as $t\to p$.

Suppose $R\ge\ln(2)/\pi$. Then $e^{\pi R}\ge 2$, and so $e^{\pi
R}-1\ge e^{\pi R}/2$.  Suppose $r>0$ and the distance from $x\in\R$
to $\{2k\pi:k\in\Z\}$ is at least $rR$, then there is $k_0\in\Z$
such that $2(k_0+1)\pi-rR\ge x\ge 2k_0\pi+rR$. Thus $2rR\le 2\pi$,
and so $r\le\pi/R$. From equation (\ref{HA-prime}), we have
$$-\HA_r'(x)=\frac{\pi^2}{r^2}\sum_{k=-\infty}^{k_0}\frac{2e^{\pi(x-2k\pi)/r}}
{(e^{\pi(x-2k\pi)/r}-1)^2}+\frac{\pi^2}{r^2}
\sum_{k=k_0+1}^{+\infty}\frac{2e^{\pi(2k\pi-x)/r}}
{(e^{\pi(2k\pi-x)/r}-1)^2}+\frac 1r$$
$$\le\frac{\pi^2}{r^2}\sum_{k=-\infty}^{k_0}\frac{2e^{\pi(2k_0\pi+rR-2k\pi)/r}}
{(e^{\pi(2k_0\pi+rR-2k\pi)/r}-1)^2}+\frac{\pi^2}{r^2}
\sum_{k=k_0+1}^{+\infty}\frac{2e^{\pi(2k\pi-(2(k_0+1)\pi-rR))/r}}
{(e^{\pi(2k\pi-(2(k_0+1)\pi-rR))/r}-1)^2}+\frac 1r$$
$$=2\frac{\pi^2}{r^2}\sum_{m=0}^\infty\frac{2e^{\pi(2m\pi+rR)/r}}
{(e^{\pi(2m\pi+rR)/r}-1)^2}+\frac
1r=\frac{\pi^2}{r^2}\sum_{m=0}^\infty\frac{4e^{\pi R+2m\pi^2/r}}
{(e^{\pi R+2m\pi^2/r}-1)^2}+\frac 1r$$
$$\le\frac{\pi^2}{r^2}\sum_{m=0}^\infty\frac{4e^{\pi R+2m\pi^2/r}}
{(e^{\pi R+2m\pi^2/r}/2)^2}+\frac 1r=16\sum_{m=0}^\infty e^{-\pi
R-2m\pi^2/r}+\frac 1r$$
$$=\frac{\pi^2}{r^2}\frac{16e^{-\pi R}}{1-e^{-2\pi^2/r}}+\frac 1r\le
\frac{\pi^2}{r^2}\frac{16e^{-\pi R}}{1-e^{-2\pi^2/(\pi/R)}}+\frac
1r\le 32\frac{\pi^2}{r^2}e^{-\pi R}+\frac 1r.$$

Let $r(t)=p_0-h(t)$ and
$R(t)=\min\{\eta(h(t))-b_-(t),b_+(t)-\eta(h(t))\}/r(t)$. Then
$r(t)\to 0$ and $R(t)\to+\infty$ uniformly in $\omega$ as $t\to p$,
and for any $x\in\til I_t$, the distance from $\eta(h(t))-x$ and
$\{2k\pi:k\in\Z\}$ is at least $r(t)R(t)$. There is $t_0\in(0,p)$
such that for $t\in[t_0,p)$, $R(t)\ge\ln(2)/\pi$ and $r(t)\le
1/(2C_0)$, where $C_0$ is as in equation (\ref{p0-h(t)}). From the
above displayed formula, we have $-H_{r(t)}'(x)\le
32\pi^2/r(t)^2e^{-\pi R(t)}+1/r$. When $t\in[t_0,p)$, $r(t)/(p-t)\ge
1/2$ by equation (\ref{p0-h(t)}), and so from the estimation of
$-\HA_{r(t)}'(x)$ and equation (\ref{estimation of mut}), we have
$$\til g_t'(\eta(h(t)))-1\le |\til\mu_t|(32\frac{\pi^2}{r(t)^2}
e^{-\pi R(t)}+\frac 1{r(t)})$$
$$\le C_0(p-t)r(t)(32\frac{\pi^2}{r(t)^2}
e^{-\pi R(t)}+\frac 1{r(t)})\le 64C_0\pi^2e^{-\pi R(t)}+C_0(p-t)\to
0,$$ as $t\to p$. Thus $\til f_t'(\xi(t))=1/\til g_t'(\eta(h(t)))\to
1$ as $t\to p$. Recall that the above argument is based on the
assumption that $T_A=p$.

Suppose
\begin{equation}\int_{p_0-h(t)}^{p-t}S_rdr\to 0\mbox{, as } t\to
p\mbox{ on the event that }T_A=p.\label{Assumption}\end{equation}
 Then $M_t\to 1$ as $t\to T_A$ on the event that $T_A=p$.
 From the Markov property, we have
$$\til\vphi_A'(0)^{5/8}\exp\left(-\frac{5}{8}
\int^{p}_{p_0}\left(S_r-\frac16\right)dr\right)=M_0=\EE[\lim_{t\to
T_A}M_t]=\PP(\{T_A=p\}).$$ Recall that $p_0$ is the modulus of
$\A_p\sem A$.  Let $K_p=\cup_{0\le t<p}K_t$. Then
\begin{equation}\PP(\{K_p\cap
A=\emptyset\})=\til\vphi_A'(0)^{5/8}\exp\left(-\frac{5}{8}
\int^{p}_{p_0}\left(S_r-\frac16\right)dr\right).\label{does not
intersect}\end{equation} If $A$ is not a smooth hull, we may find a
sequence of smooth hulls $A_n$ that approaches $A$. Then
$\til\vphi_{A_n}'(0)\to\til\vphi_A'(0)$ and the modulus of $\A_p\sem
A_n$ tends to the modulus of $\A_p\sem A$, so equation (\ref{does
not intersect}) still holds.

Now suppose $B$ is a hull in $\A_{p_0}$ w.r.t.\ $\CC_{p_0}$. Let
$D=A\cup\vphi_A^{-1}(B)$. Then $D$ is a hull in $\A_p$ w.r.t.\
$\CC_p$. Let $p_1$ be the modulus of $\A_p\sem D$, which is also the
modulus of $\A_{p_0}\sem B$. Then $\vphi_D=\vphi_B\circ\vphi_A$, so
$\til\vphi_D=\til\vphi_B\circ\til\vphi_A$ and
$\til\vphi_D'(0)=\til\vphi_B'(\til\vphi_A(0))\til\vphi_A'(0)=\til\vphi_B'(0)
\til\vphi_A'(0)$. From the last paragraph,
$$\PP(K_p\cap D=\emptyset\})=
\til\vphi_D'(0)^{5/8}\exp\left(-\frac{5}{8}
\int^{p}_{p_1}\left(S_r-\frac16\right)dr\right).$$ Thus
$$\PP(\{K_p\cap D=\emptyset\}|
\{K_p\cap A=\emptyset\})=\til\vphi_B'(0)^{5/8}\exp\left(-\frac{5}{8}
\int^{p_0}_{p_1}\left(S_r-\frac16\right)dr\right).$$ If $L_s$, $0\le
s<p_0$, are standard annulus SLE$_{8/3}$ hulls of modulus $p_0$, and
$L_{p_0}=\cup_{0\le s<p_0}L_s$, then
$$\PP(\{L_{p_0}\cap B=\emptyset\})=\til\vphi_B'(0)^{5/8}\exp\left(-\frac{5}{8}
\int^{p_0}_{p_1}\left(S_r-\frac16\right)dr\right).$$
$$=\PP(\{K_p\cap D=\emptyset\}|
\{K_p\cap A=\emptyset\})=\PP(\{\vphi_A(K_p)\cap
B=\emptyset\}|\{K_p\cap A=\emptyset\}).$$ Thus conditioned on the
event that $K_p\cap A=\emptyset$, $\vphi_A(K_p)$ has the same
distribution as $L_{p_0}$. Then we proved the restriction property
of annulus SLE$_{8/3}$ under the assumption (\ref{Assumption}).

Unfortunately, the assumption (\ref{Assumption}) is actually always
false. From equation (\ref{change HA}) one may compute that $S_r$ is
of order $\Theta(1/r^2)$ as $r\to 0$. From (\ref{p0-h(t)}),
$(p-t)-(p_0-h(t))=|\mu_t|$ is of order $O((p-t)^2)$. In fact, one
could prove that it is of order $\Theta((p-t)^2)$. So
$\int_{p_0-h(t)}^{p-t} S_rdr$ is uniformly bounded away from $0$.
Thus it does not tend to $0$ as $t\to p$. Therefore we guess that
annulus SLE$_{8/3}$ does not satisfy the restriction property.

\vskip 3mm

Recently, Robert O. Bauer  studied in \cite{recent} a process
defined in a doubly connected domain obtained  by conditioning a
chordal SLE$_{8/3}$ in a simply connected domain to avoid an
interior contractible compact subset. The process describes a random
simple curve connecting two prime ends of a doubly connected domain
that lie on the same side, so it is different from the process we
study here. That process automatically satisfies the restriction
property from the restriction property of chordal SLE$_{8/3}$. And
it satisfies conformal invariance because the set of boundary hulls
generates the same $\sigma$-algebra as the Hausdorff metric on the
space of simple curves.

\vskip 4mm

\no{\bf Acknowledgement:} I would like to thank my PhD advisor
Nikolai Makarov for his instruction on this work.

\end{document}